\documentclass{amsart}

\usepackage{amssymb}
\usepackage{amsmath}
\usepackage{bbm}
\usepackage{graphicx}

\usepackage{amsfonts}

\newcommand{\assign}{:=}
\newcommand{\nonesep}{}
\newcommand{\tmdummy}{$\mbox{}$}
\newcommand{\tmem}[1]{{\em #1\/}}

\newcommand{\tmop}[1]{\ensuremath{\operatorname{#1}}}
\newenvironment{itemizedot}{\begin{itemize}
    }{\end{itemize}}

\def\rb{\mathbf{r}} 
\def\In{\tmop{In}} 
\def\MIS{\tmop{\mathrm{\textsc{mis}}}} 
 
\def\Tc{\mathcal{T}} 
\def\NN{\mathbbm{N}}
\def\RR{\mathbbm{R}}

\newtheorem{theorem}{Theorem}[section]
\newtheorem{lemma}[theorem]{Lemma}
\newtheorem{proposition}[theorem]{Proposition}
\newtheorem{corollary}[theorem]{Corollary}
\newtheorem{algorithm}[theorem]{Algorithm}
\theoremstyle{definition}
\newtheorem{definition}[theorem]{Definition}
\newtheorem{example}[theorem]{Example}

\theoremstyle{remark}
\newtheorem{remark}[theorem]{Remark}

\numberwithin{equation}{section}
\begin{document}

\title{\scshape On the dimension of spline spaces on planar T-meshes}
\author[Bernard Mourrain]{Bernard Mourrain,\\
GALAAD, INRIA M\'editerran\'ee\\
\texttt{Bernard.Mourrain@inria.fr}}

\begin{abstract}
  We analyze the space $\mathcal{S}_{m, m'}^{\rb} (\Tc)$ of bivariate
  functions that are piecewise polynomial of bi-degree $\leqslant (m,
  m')$ and of smoothness $\rb$ along the interior edges of a planar
  T-mesh $\Tc$.  We give new combinatorial lower and upper bounds for
  the dimension of this space by exploiting homological techniques. We
  relate this dimension to the weight of the maximal interior segments
  of the T-mesh, defined for an ordering of these maximal interior segments.
  We show that the lower and upper bounds coincide, for high enough
  degrees or for hierarchical T-meshes which are enough regular.  We
  give a rule of subdivision to construct hierarchical T-meshes for
  which these lower and upper bounds coincide.
  Finally, we illustrate these results by analyzing spline spaces of
  small degrees and smoothness.
\end{abstract}

\maketitle

\section*{Introduction}
Standard parametrisations of surfaces in Computer Aided Geometric Design are
based on tensor product B-spline functions, defined from a grid of nodes over
a rectangular domain \cite{farin:book}.  These representations are easy to control but their
refinement has some drawback.  Inserting a node in one direction of the
parameter domain implies the insertion of several control points in the other
directions. If for instance, regions along the diagonal of the parameter
domain should be refined, this would create a fine grid in some regions where it is
not needed. To avoid this problem, while extending the standard tensor
product representation of CAGD, spline functions associated to a subdivision with
T-junctions instead of a grid, have been studied. Such a T-mesh
is a partition of a domain $\Omega$ into
axis-aligned boxes, called the cells of the T-mesh.

The first type of T-splines introduced in
\cite{Sederberg03,Sederberg04}, are defined by blending functions
which are products of univariate B-spline basis functions associated to
some nodes of the subdivision.  They are piecewise polynomial
functions, but the pieces where these functions are polynomial do not
match with the cells of the T-subdivision.  Moreover, there is no
proof that these piecewise polynomial functions are linearly
independent. Indeed, \cite{BuChSa10} shows that in some cases, these
blending T-spline functions are not linearly independent.  Another
issue related to this construction is that there is no
characterization of the vector space spanned by these functions.  For
this reason, the partition of unity property which is useful in CAGD is not
available directly in this space. The spline functions have to be
replaced by piecewise rational functions, so that these piecewise
rational functions sum to $1$. However, this construction complexifies
the practical use of such T-splines.

Being able to describe a basis of the vector space of piecewise
polynomials of a given smoothness on a T-mesh is an important
but non-trivial issue. It yields a construction of piecewise
polynomial functions on the T-subdivision which form a partition of
unity so that the use of piecewise rational functions is not
required. It has also a direct impact in approximation problems such as
surface reconstruction \cite{deBoor01} or isogeometric analysis
\cite{hughes:CMAME2005}, where controlling the space of functions used
to approximate a solution is critical. In CAGD, it also provides more
degrees of freedom to control a shape.  This explains why further
works have been developed to understand better the space of piecewise
polynomial functions with given smoothness on a T-subdivision.

To tackle these issues special families of splines on
T-meshes have been studied. 
In {\cite{deng06}}, {\cite{deng08}}, 
these splines are piecewise polynomial functions on a hierarchical
T-subdivision. They are called PHT-splines (Polynomial Hierarchical
T-splines). Dimension formulae of the spline space on such a
subdivision have been proposed when the degree is high enough compared
to the smoothness {\cite{deng06}, \cite{HuDeFeChe06},
  \cite{LiWaZha06}} and in some cases for biquadratic $C^1$ piecewise
polynomial functions {\cite{Deng08TSPL22}}. The construction of a
basis is described for bicubic $C^1$ spline spaces in terms of the
coefficients of the polynomials in the Bernstein basis attached to a
cell. When a cell is subdivided into 4 subcells, the Bernstein
coefficients of the basis functions of the old level are modified and
new linearly independent functions are introduced, using Bernstein
bases on the cells at the new level.

In this paper, we analyse the dimension of the space $\mathcal{S}_{m, m'}^{\rb}(\Tc)$ 
of bivariate functions that are piecewise polynomial of
bi-degree $\leqslant (m, m')$ of smoothness $\rb$ along the interior
edges of a general planar T-mesh $\Tc$, where $\rb$ is a smoothness distribution on
$\Tc$. 

As we will see, computing this dimension reduces to compute the dimension of the kernel of a certain
linear map (namely the map $\tilde{\partial}_{2}$ introduced in Section
\ref{sec:2}).

Thus for a given T-mesh, a given smoothness distribution $\rb$ and
a given bi-degree $(m,m')$, it is possible to compute the dimension of
$\mathcal{S}_{m, m'}^{\rb} (\Tc)$ by linear algebra tools (see eg. a
software implementation developed by 
P. Alfed\footnote{\texttt{http://www.math.utah.edu/$\sim$pa/MDS/index.html}} 
for such computations).
We would like to avoid a case-by-case treatment and to
describe this dimension in terms of combinatorial quantities attached to $\Tc$ 
and easy to evaluate. As shown in \cite{Li:2011:IDS:2038067.2038142} 
or \cite{Berdinsky:2012}, the dimension may also depend on the
geometry of the T-mesh and not just on its topology. 
This explains why it is not always possible to provide a purely
combinatorial formula for the dimension of $\mathcal{S}_{m, m'}^{\rb} (\Tc)$. 

The main results in this paper are 
\begin{itemize}
 \item a description of the dimension in terms of a {\em combinatorial}
   part that depends only on the topology of the T-mesh 
and an {\em homological} part that takes into
   account the fact that the dimension may depend on the
   geometry of the T-mesh (Theorem \ref{thm:dim}); 

 \item {\em combinatorial} lower and upper bounds on the
   dimension that are easy to evaluate (Theorem \ref{thm:dim:bound});

 \item sufficient conditions under which the lower and upper bounds coincide 
so that the dimension depends only on the topology of the T-mesh 
(Theorem \ref{thm:weighted}). 
\end{itemize}
 
We proceed as follows. By extending homological techniques developed
in \cite{b-htssg-88} and \cite{Schenck1997535}, we obtain
combinatorial lower and upper bounds on the dimension of these spline
spaces for general T-subdivisions.  We relate the upper bound to the
maximal interior segments and their weights and show that the lower
and upper bounds coincide for $T$-meshes which are enough regular.
Namely, if a T-mesh is $(m+1, m'+1)$-weighted, the dimension depends
directly on the number of faces, interior edges and interior points.
In particular, we obtain the dimension formula for a constant
smoothness distribution $\rb= (r,r')$ with $m \geq 2 r + 1$ and $m'
\geq 2 r' + 1$, providing a new proof of a result also available in
{\cite{deng06}, \cite{HuDeFeChe06}, \cite{LiWaZha06}} for a
hierarchical T-mesh. The algebraic approach gives an homological
interpretation of the method called Smoothing Cofactor-Conformality
method in \cite{Wang01}. It allows us to generalize the dimension
formulae obtained by this technique \cite{LiWaZha06},
\cite{HuDeFeChe06}.  We also give a rule of subdivision to construct
hierarchical T-meshes for which the lower and upper bounds
coincide. As a consequence, we can recover the dimension of the space
of Locally Refined splines described in {\cite{DoLyPe10}}.  We do not
consider the problem of constructing explicit bases for these spline
spaces, which will be analyzed separately.

In the first section, we recall the notations and the polynomial properties
which are needed in the following. 
Section \ref{sec:2} describes the chain complex associated to the
spline space and analyzes its homology.
In Section \ref{sec:3}, we give lower and upper bounds on the dimension of the spline
space and analyze cases where these bounds are coincide.
Section \ref{sec:4} deals with the properties of hierarchical
$T$-meshes, obtained by recursive subdivisions of cells.
In the last section, we analyse some examples for small degree and
smoothness.

\section{Planar T-splines}\label{sec:1}

In the following, we will deal with notions which are of topological
and algebraic nature. We start by the topological definitions.

\subsection{T-meshes}
For any set $S\subset \RR^{2}$, $\overline{S}$ is its closure for the
usual topology, $S^{\circ}$ its relative interior, $\partial S$ its boundary.

We define a T-mesh $\Tc$ of $\RR^{2}$ as:
\begin{itemize}
 \item a finite set $\Tc_{2}$ of closed axis-aligned rectangles of
$\RR^{2}$,
 \item a finite set $\Tc_{1}$ of closed axis-aligned segments
   included in $\cup_{\sigma\in \Tc_{2}} \partial \sigma$,
 \item a finite set of points $\Tc_{0} = \cup_{\tau\in \Tc_{1}} \partial \tau$,
\end{itemize}
  such that
\begin{itemize}
\item For $\sigma\in \Tc_{2}$, $\partial \sigma$ is the finite union of elements of $\Tc_{1}$.
\item For $\sigma,\sigma'\in \Tc_{2}$ with $\sigma\neq \sigma'$, 
$\sigma\cap \sigma'=\partial \sigma\cap \partial \sigma'$ is the finite union of
 elements of $\Tc_{1}\cup \mathcal{T}_{0}$.
 \item For $\tau,\tau'\in \Tc_{1}$ with $\tau\neq \tau'$, $\tau\cap
   \tau' =\partial \tau \cap \partial\tau' \subset \Tc_{0}$.
\end{itemize}
We denote by $\Omega =\cup_{\sigma\in \Tc_{2}} \sigma \subset \RR^2$
and call it the domain of the T-mesh $\Tc$.  

The elements of $\Tc_{2}$ are called $2$-faces or cells and their number is denoted $f_2$.

The elements of $\Tc_{1}$  are called $1$-faces or edges.
An element of $\Tc_{1}$ is called an interior edge if it intersects
$\Omega^{\circ}$.
It is called a boundary edge otherwise.
The set of interior edges is denoted by $\Tc_1^o$. The number of edges
in $\Tc_1$ is $f_1$ and the number of interior edges is $f_{1}^{o}$.

An edge parallel to the first (resp. second) axis of $\RR^{2}$ is called horizontal
(resp. vertical). Let $\Tc_{1}^{o,h}$ (resp. $\Tc_{1}^{o,v}$) be the set of horizontal
(resp. vertical) interior edges and $f_1^h$ (resp. $f_1^v$) the number of interior
horizontal (resp. vertical) edges. Then, the number of interior edges is
$f_1^o = f_1^h + f_1^v$.
 
The elements of $\Tc_{0}$ are called $0$-faces or vertices.
A vertex is interior if it is in $\Omega^{\circ}$. It is a boundary
vertex otherwise. 
The set of interior vertices is denoted $\Tc_0^o$.
We denote by $f_0$ be the number of vertices of $\Tc_0$ and by $f_0^o$ be the
number of interior vertices.

A vertex is a crossing vertex if it is an interior vertex and belongs
to $4$ distinct elements of $\Tc_{1}$.  A vertex is a T-vertex if it
is an interior vertex and belongs to exactly $3$ distinct elements of
$\Tc_{1}$.  Let $f_{0}^{+}$ (resp. $f_0^T$) be the number of crossing
(resp. T) vertices.  A boundary vertex is a vertex in
$\Tc_{0}\cap \partial \Omega$. The number of boundary vertices is
$f_{0}^{b}$.  A vertex is a corner vertex if it belongs to
$\partial \Omega$ and to a vertical and a horizontal boundary edge.

To simplify the definitions and remove redundant edges, 
we will assume that {\em any vertex $\gamma\in \Tc_{0}$ belongs at
  least to one horizontal edge $\tau_{h} (\gamma) \in \Tc_{1}$ and one vertical edge
  $\tau_{v} (\gamma) \in \Tc_{1}$}.

We denote by $\nu_{h} (\Tc) =\{s_{1}, \ldots, s_{l}\}\subset \RR$ 
(resp. $\nu_{v}(\Tc) =\{t_{1}, \ldots, t_{m}\} \subset \RR$)
the set of first (resp. second) coordinates of
the points in vertical (horizontal) segments $\in \Tc_{1}^{v}$
(resp. $\in \Tc_{1}^{h}$).
The elements of $\nu_{h} (\Tc)$ (resp. $\nu_{v} (\Tc)$) are called the horizontal
(resp. vertical) nodes of the T-mesh $\Tc$.

\begin{example}\label{ex:1}
Let us illustrate the previous definitions on the following T-mesh:
\begin{center}
\begin{center}
    \includegraphics[height=4.5cm]{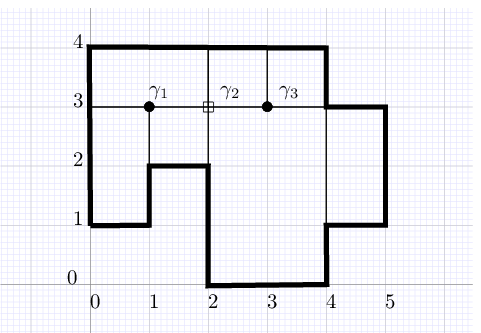}
\end{center}
\end{center}
In this example, there are $f_{2}= 7$ rectangles, $f_{1}^{o}=9$ interior
edges such that $f_{1}^{h}=4$ are horizontal
and $f_{1}^{v}=5$ are vertical.
There are $f_{0}^{o}=3$ interior points $\gamma_{1},
\gamma_{2},\gamma_{3}$; $\gamma_{1}, \gamma_{3}$ are T-vertices and
$\gamma_{2}$ is a crossing vertex. 
There are $f_{0}^{b}=15$ boundary vertices and $12$ corner vertices.

The horizontal nodes are $\nu_{h} (\Tc) =\{0,
\ldots, 5\}$ and the vertical nodes are $\nu_{v} (\Tc) =\{0, \ldots, 4\}$.
\end{example}

For our analysis of spline spaces on T-meshes, we assume the following:

\noindent{}\textbf{Assumption:} {\em The domain $\Omega$ is simply
connected and $\Omega^{o}$ is connected.}

This implies that $\Omega$ has one connected component with no ``hole''
and that the number of boundary edges through a boundary vertex is $2$.
 
\subsection{T-splines}
We are going now to define the spaces of piecewise polynomial
functions on a T-mesh, with bounded degrees and given smoothness. An element
in such a space is called a spline function.

\begin{definition}
A {\em smoothness distribution} on a T-mesh $\Tc$ is a pair of maps
$(\rb_{h},\rb_{v})$ from $(\nu_{h} (\Tc) \times \nu_{v} (\Tc))$ to $\NN\times \NN$.
\end{definition}

For convenience, we will define the smoothness map $\rb$ on $\Tc_{1}$ as follows: for any $\tau \in
\Tc_{1}^{v}$ (resp. $\tau \in \Tc_{1}^{h}$), 
$\rb (\tau)= \rb_{h} (s)$ 
(resp. $\rb (\tau)= \rb_{v} (t)$)
where $s \in \nu_{h} (\Tc)$ 
(resp. $t\in \nu_{v} (\Tc)$)
is the first
(resp. second) coordinate of any point of $\tau$.
We will also define the horizontal and vertical smoothness on
$\Tc_{0}$ as follows: 
for any $\gamma = (s,t) \in
\Tc_{0}$, $\rb_{h} (\gamma) = \rb_{h} (s)$ and 
$\rb_{v} (\gamma) = \rb_{v} (t)$.

For $r,r'\in \NN$, we say that $\rb$ is a constant smoothness
distribution equal to $(r,r')$ if 
$\forall s \in \nu_{h} (\Tc), \rb (s)=r$,
$\forall t \in \nu_{v} (\Tc), \rb (t)=r'$. 

Let $R=\RR[s,t]$ be the polynomial ring with coefficients in $\RR$.
For $m,m'\in \NN$, we denote by $R_{m,m'}=R_{(m,m')}$ the vector of polynomials in $R$
of degree $\le m$ in $s$ and $\le m'$ in $t$. An element of $R_{m,m'}$
is of bi-degree $\leqslant (m,m')$.


The goal of this paper is to analyse the dimension of the 
space of splines of bi-degree $\le (m,m')$ and of smoothness $\rb$ on
a T-mesh $\Tc$, that we define now.
\begin{definition}
Let $\Tc$ be a T-mesh and $\rb$ a node smoothness distribution.
We denote by $\mathcal{S}_{m, m'}^{\rb} (\Tc)$ the vector space of
  functions which are polynomials in $R_{m, m'}$ on each cell $\sigma \in
  \Tc_2$ and of class $\rb (\tau)$ in $s$ (resp. in $t$) at any point 
of $\tau \cap \Omega^{o}$ 
if $\tau$ is a vertical (resp. horizontal) interior edge.
\end{definition}

We will say that $f \in \mathcal{S}_{m, m'}^{\rb} (\Tc)$ is of (continuity) class $C^{\rb}$
on $\Tc$. 
We notice that the boundary edges and their smoothness are not involved in the
characterization of a spline function.

\begin{example} We consider again the T-mesh of Example \ref{ex:1}.
If we take the node smoothness distribution 
$\rb_{h} (1)= 1$, $\rb_{h} (2)= 0$, $\rb_{h} (3)= \rb_{h}
(4)= \rb_{h} (5)= 1$,  and $\rb_{v}$ constant equal to $1$, then 
$\mathcal{S}_{3,3}^{\rb} (\Tc)$ is the vector space of bicubic piecewise
polynomial functions on $\Tc$ which are $C^{1}$ in $s$ for $s<2$ and $s>2$,
continuous for $s=2$ and $C^{1}$ in $t$ in $\Omega^{o}$.
\end{example}
 
\subsection{Polynomial properties}

We recall here basic results on the dimension of the vector spaces involved in
the analysis of \ $\mathcal{S}_{m, m'}^{\rb} (\Tc)$:

For any $\tau \in \Tc_{1}$, let $l_{\tau} \in R$ be a
non-zero polynomial (of degree $1$) defining the line supporting the edge $\tau$. 
Let $\Delta^{\rb} (\tau) = l_{\tau}^{\rb (\tau) + 1}$. 
We denote by $\mathfrak{I}^{\rb} (\tau) =
(\Delta^{\rb} (\tau))$ the ideal generated by the polynomial
$\Delta^{\rb} (\tau)\in R$ and by $\mathfrak{I}^{\rb}_{m, m'} (\tau)
=\mathfrak{I}^{\rb} (\tau) \cap R_{m, m'}$ its part of bi-degree $\leq (m,m')$.
Notice that $\mathfrak{I}^{\rb}(\tau)$ defines the line supporting the
egde $\tau$ with multiplicity $(\rb (\tau) + 1)$.
By definition, two horizontal (resp. vertical) edges $\tau_1, \tau_2$ which
share a point define the same ideal $\mathfrak{I}^{\rb} (\tau_1)
=\mathfrak{I}^{\rb} (\tau_2)$.
We define the bi-degree $\delta$ for any egde $\tau \in \Tc_{1}$ as follows:
\begin{itemize}
\item $\delta (\tau) = (\rb (\tau) + 1, 0)$ if $\tau$ is vertical, 
\item $\delta (\tau) = (0, \rb (\tau) + 1)$ if $\tau$ is horizontal.
\end{itemize}

Let $\mathfrak{I}^{\rb} (\gamma) =\mathfrak{I}^{\rb}
 (\tau_v) +\mathfrak{I}^{\rb} (\tau_h) = (\Delta^{\rb} (\tau_v),
 \Delta^{\rb} (\tau_h))$
where $\tau_{v}, \tau_{h}\in \Tc_{1}$ are vertical and horizontal
edges such that $\tau_{v}\cap \tau_{h}= \{\gamma\}$.
The ideal $\mathfrak{I}^{\rb}(\gamma)$ defines the point $\gamma$ with
multiplicity $(\rb_{h} (\gamma)  + 1) \times (\rb_{v} (\gamma) + 1)$.
We denote by
 $\mathfrak{I}^{\rb}_{m, m'} (\gamma) =\mathfrak{I}^{\rb}_{m, m'}
 (\tau_v) +\mathfrak{I}^{\rb}_{m, m'} (\tau_h)$. Notice that these
 definitions are independent of the choice of the vertical
 edge $\tau_v$ and horizontal edge $\tau_h$ which contain $\gamma$.
The bi-degree of a vertex $\gamma\in \Tc_{0}$ is $\delta (\gamma)=
(\rb_{h} (\gamma)  + 1, \rb_{v} (\gamma) + 1)$.

Here are the basic dimension relations that we will use to analyse the
spline functions on a T-mesh.
\begin{lemma}
  \label{lem:dim}{\tmdummy}
  
  \begin{itemizedot}
    \item $\dim R_{m, m'} = (m + 1) \times (m' + 1)$.
    
    \item $\dim \left( R_{m, m'} /\mathfrak{I}^{\rb}_{m, m'} (\tau)
      \right) = \left\{ 
    \begin{array}{ll}
      (m + 1) \times (\min (\rb(\tau),m') + 1) & \tmop{if} \tau \in \Tc_{1}^{h}\\
      (\min (\rb(\tau),m) + 1) \times (m' + 1) & \tmop{if} \tau \in \Tc_{1}^{v}
    \end{array}  \right.$
    
    \item $\dim \left( R_{m, m'} /\mathfrak{I}^{\rb}_{m, m'}
        (\gamma) \right) = (\min (\rb_{h}(\gamma),m) + 1)
    \times (\min (\rb_{v} (\gamma),m') + 1)$ for all $\gamma\in \Tc_{0}$.
  \end{itemizedot}
\end{lemma}
\begin{proof} To obtain these formulae, we directly check that 
\begin{itemize}
 \item a basis of $R_{m,m'}$ is the set of monomials $s^{i} t^{j}$
   with 
$0 \le i \le m$, $0 \le j\le m'$;
 \item a basis of $R_{m,m'}/\mathfrak{I}^{\rb}_{m, m'} (\tau)$ is the
   set of monomials $s^{i} t^{j}$ with
$0 \le i \le m$ and $0 \le j\le \min (\rb (\tau),m')$ if $\tau \in \Tc_{1}^{h}$
(resp. $0 \le i \le \min (\rb (\tau),m)$, $0 \le j\le m'$ if $\tau \in \Tc_{1}^{v}$);
 \item a basis of $R_{m,m'}/\mathfrak{I}^{\rb}_{m, m'} (\gamma)$ is
   the set of monomials $s^{i} t^{j}$ with 
$0 \le i \le \min (\rb_{h} (\gamma),m)$, $0 \le j\le \min (\rb_{v}
(\gamma),m')$.
\end{itemize} 
since the ideal of an edge $\tau \in \Tc_{1}$ is up to a translation
$(s^{\rb(\tau)+1})$  or $(t^{\rb(\tau)+1})$.
\end{proof}
 
An algebraic way to characterise the $C^{\rb}$-smoothness is given by the
next lemma:
\begin{lemma}[\cite{b-htssg-88}]
  \label{lem:reg}Let $\tau \in \Tc_1$ and let $p_1, p_2$ be two
  polynomials. Their derivatives coincide on $\tau$ up to order $\rb (\tau)$ iff
  $p_1 - p_2 \in \mathfrak{I}^{\rb} (\tau)$.
\end{lemma}

In the following, we will need algebraic properties on univariate polynomials.
We denote by $U=\RR[u]$ the space of univariate polynomials in the variable $u$ with
coefficients in $\RR$.
Let $U_{n}$ denote the space of polynomials of $U$ of degree $\leq n$.
For a polynomial $g\in U$ of degree $d$ and an integer $n\geq d$, $g\, U_{n-d}$ is the 
vector space of multiples of $g$ which are of degree $\leq n$.
For polynomials $g_{1}, \ldots, g_{k}\in U$ respectively of degree
$d_{1},\ldots, d_{k}$ and an integer $n\geq \max_{i=1,\ldots,k}
d_{i}$, 
$\sum_{i=1}^{k}\, g_{i}\, U_{n-d_{i}}$ is the vector space of sums of multiples of $g_{i}$ of degree $\leq n$.

We will use the apolar product defined on $U_{n}$ by:
\[ \langle f, g \rangle_n = \sum_{i = 0}^n  \binom{n}{i} f_i g_i \]
where $f = \sum_{i = 0}^n f_i u^i $, $g = \sum_{i = 0}^n g_i u^i \in
\RR[u]_n$. One of the properties that we will need is the following
{\cite{Salmon1885}, \cite{KungRota84}, \cite{EhrRota93}}:
\begin{lemma}
  \label{lem::polar}Let $g \in U_n$, $d < n$ and $a \in
  \RR$. Then $g$ is orthogonal to $(u - a)^d U_{n - d}$ for the apolar
  product iff
 \[ \partial^k g (a) = 0, k = 0, \ldots, n - d. \]
\end{lemma}

\begin{proposition}
  \label{dim:apolar}Let $a_1, \ldots, a_l \nonesep \in \RR$ be $l$
  distinct points and $d_{1}, \ldots, d_{l} \in \NN$ . Then
  \[ \dim \left(  \sum_{i = 1}^l (u - a_i)^{d_{i}} U_{n - d_{i}} \right) =
     \min (n + 1, \sum_{i=1}^{l} n - d_{i} + 1) . \]
\end{proposition}
\begin{proof}
  In order to compute the dimension of $V \assign \sum_{i = 1}^l (u - a_i)^d
  U_{n - d} \subset U_n$, we compute the dimension
  of the orthogonal $V^{\bot}$ in $U_n$ of $V$ for the apolar
  product. Let $g \in U_n$ be an element of the orthogonal
  $V^{\bot}$ of $V$. By lemma \ref{lem::polar}, $\partial^k g (a_i) = 0, k = 0,
  \ldots, n - d_{i}, i = 1, \ldots, l$. In other words, $g$ is divisible by $(u -
  a_i)^{n - d_{i} + 1}$ for $i = 1, \ldots, l$. As the points $a_i$ are distinct,
  $g$ is divisible by
  \[ \Pi \assign \prod_{i = 1}^l (u - a_i)^{n - d_{i} + 1} . \]
  Conversely, any multiple of $\Pi$ of degree $\leqslant n$ is in
  $V^{\bot}$. Thus $V^{\bot} = (\Pi)\cap U_n$.
  
  This vector space $V^{\bot}$ of multiples of $\Pi$ in degree $n$ is of
  dimension $\max (0, n + 1 - \deg (\Pi))$, so $V$ is of dimension
  \[ n + 1 - \max (0, n + 1 - \deg (\Pi)) = \min (n + 1, \deg (\Pi)) =
     \min (n + 1, \sum_{i=1}^{l} n - d_{i} + 1). \]
\end{proof}
We are going to use an equivalent formulation of this result:
\begin{equation}
  \dim \left( U_n / \sum_{i = 1}^l (u - a_i)^d_{i} U_{n - d_{i}} \right)
  = (n + 1 - \sum_{i=1}^{l} n - d_{i} + 1)_+
\end{equation}
where for any $a \in \mathbb{Z}$, $a_{+}=\max (0, a)$.
A similar result is proved in \cite{LiWaZha06}[Lemma 2] when all $d_{i}$ are
equal to $d$, by analyzing the coefficient matrix of generators
of $\sum_{i = 1}^l (u - a_i)^d U_{n - d}$. 

\subsection{Maximal interior segments}

In order to simplify the analysis of $\mathcal{S}_{m, m'}^{\rb}
(\Tc)$, we introduce the following definitions:

For any interior edge $\tau \in \Tc_1^o$, we define $\rho(\tau)$ as the
 {\tmem{maximal segment}} made of edges $\in \Tc_1^o$
 of the same direction as $\tau$, which contains $\tau$ and such that their
 union is connected. We say that the maximal segment $\rho(\tau)$ is interior if it does not
 intersect the boundary of $\Omega$.
 
As all the edges belonging to a maximal segment $\rho$ have the same
 supporting line, we can define $\Delta^{\rb} (\rho) = \Delta^{\rb} (\tau)$
 for any edge $\tau$ belonging to $\rho= \rho(\tau)$.
 
The set of all maximal interior segments is denoted
by $\MIS (\Tc)$. The set of horizontal (resp. vertical) maximal
 interior segments of $\Tc$ is denoted by $\MIS_h
 (\Tc)$ (resp. $\MIS_v (\Tc)$).
 
The degree of $\rho \in \MIS (\Tc) \overline{}$ is by definition $\delta
 ( \rho) = \delta (\tau)$ for any $\tau \subset \rho$.

For each interior vertex $\gamma \in \Tc_0^o$, which is the
 intersection of an horizontal edge $\tau_h \in \Tc_1^o$ with a
 vertical edge $\tau_v \in \Tc_1^o$, let $\rho_h (\gamma)$ (resp. $\rho_v
 (\gamma)$) is the corresponding horizontal (resp. vertical) maximal
 segment. We denote by $\Delta_h^{\rb} (\gamma)$
 (resp. $\Delta_v^{\rb} (\gamma)$) the equations of the
 corresponding supporting lines to the power $\rb_{h} (\gamma) + 1$ (resp. $\rb_{v} (\gamma) + 1$).

Notice that the intersection of two distinct maximal interior segments is
either a T-vertex or a crossing vertex.

We say that $\rho \in \MIS (\Tc)$ is {\tmem{blocking}} $\rho' \in
\MIS (\Tc)$ if one of the end points of $\rho'$ is in the interior
of $\rho$.

\begin{example}\label{ex:1.9}
 In the figure below, the maximal interior edges are indicated by plain
 segments.
  \begin{center}
    \includegraphics{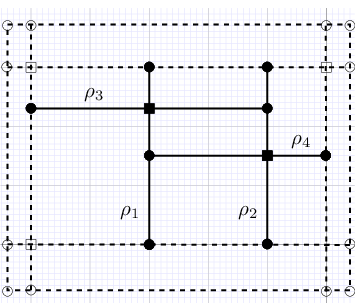}
  \end{center}
In this example, $\rho_1$ is blocking $\rho_4$ and $\rho_2$ is blocking $\rho_3$.
\end{example}

\section{Topological chain complexes}\label{sec:2}
In this section, we describe the tools from algebraic topology, that we
will use. For more details, see eg. {\cite{Spanier66}}, {\cite{Hatcher02}}.

\subsection{Definitions}

We consider the following complexes:
{\small
\[ 
\!\!\!\!\!\!\!\!\!\!\!\!\!\!\!\!\!\!\!\!\!\!\!\!
\begin{array}{lllllllll}
     &  &  &  & 0 &  & 0 &  & \\
     &  &  &  & \downarrow &  & \downarrow &  & \\
     \mathfrak{I}_{m, m'}^{\rb} (\Tc^o) : &  & 0 & \rightarrow &
     \bigoplus_{\tau \in \Tc_1^o} [\tau]\mathfrak{I}^{\rb}_{m, m'}
     (\tau) & \rightarrow & \bigoplus_{\gamma \in \Tc_0^o}
     [\gamma]\mathfrak{I}^{\rb}_{m, m'} (\gamma) & \rightarrow & 0\\
     &  & \downarrow &  & \downarrow &  & \downarrow &  & \\
     \mathfrak{R}_{m, m'} (\Tc^o) : &  & \bigoplus_{\sigma \in
     \Tc_2} [\sigma] R_{m, m'} & \rightarrow & \bigoplus_{\tau \in
     \Tc_1^o} [\tau] R_{m, m'} & \rightarrow & \bigoplus_{\gamma \in
     \Tc_0^o} [\gamma] R_{m, m'} & \rightarrow & 0\\
     &  & \downarrow &  & \downarrow &  & \downarrow &  & \\
     \mathfrak{S}_{m, m'}^{\rb} (\Tc^o) : &  & \bigoplus_{\sigma \in
     \Tc_2} [\sigma] R_{m, m'} & \rightarrow & \bigoplus_{\tau \in
     \Tc_1^o} [\tau] R_{m, m'} /\mathfrak{I}^{\rb}_{m, m'} (\tau) &
     \rightarrow & \bigoplus_{\gamma \in \Tc_0^o} [\gamma] R_{m, m'}
     /\mathfrak{I}^{\rb}_{m, m'} (\gamma) & \rightarrow & 0\\
     &  & \downarrow &  & \downarrow &  & \downarrow &  & \\
     &  & 0 &  & 0 &  & 0 &  & 
   \end{array} \]
}
The different vector spaces of these complexes are obtained as the components
in bi-degree $\leqslant (m, m')$ of $R$-modules generated by (formal)
independent elements $[\sigma], [\tau], [\gamma]$ indexed respectively by the
faces, the interior edges and interior points of $\Tc$. An oriented
edge $\tau \in \Tc_1$ is represented as: $[\tau] = [a b]$ where $a,
b \in \Tc_0$ are the end points. The opposite edge is represented by
$[b a]$. By convention, $[b a] = - [a b]$.

The maps of the complex $\mathfrak{R}_{m, m'} (\Tc^o)$ are defined as
follows:
\begin{itemizedot}
  \item for each face $\sigma \in \Tc_2$ with its counter-clockwise
  boundary formed by edges $\tau_1 = a_1 a_2, \ldots, \tau_l = a_l a_1$,
  $\partial_2 (\sigma) = [\tau_1] \oplus \cdots \oplus [\tau_l] = [a_1 a_2]
  \oplus \cdots \oplus [a_l a_1]$;
  
  \item for each interior edge $\tau \in \Tc_1^o$ from $\gamma_1$ to $\gamma_2 \in \Tc_0$, $\partial_1 ([\tau]) = [\gamma_2] -
  [\gamma_1]$ where $[\gamma] = 0$ if $\gamma \not\in \Tc_0^o$;
  
  \item for each interior point $\gamma \in \Tc_{0}^{o}$, $\partial_0 ([\gamma]) =
  0$.
\end{itemizedot}
By construction, we have $\partial_i \circ \partial_{i + 1} = 0$ for $i = 0,
1$. The maps of the complex $\mathfrak{I}_{m, m'}^{\rb} (\Tc^o)$ are
obtained from those of $\mathfrak{R}_{m, m'} (\Tc^o)$ by restriction,
those of the complex $\mathfrak{S}_{m, m'}^{\rb} (\Tc^o)$ are
obtained by taking the quotient by the corresponding vector spaces of
$\mathfrak{I}_{m, m'}^{\rb} (\Tc^o)$. The corresponding
differentials of the complex are denoted $\bar{\partial}_i$.

For each column of the diagram, the vertical maps are respectively the
inclusion map and the quotient map.

The complex $\mathfrak{R}_{m, m'} (\Tc^o)$ is also known as the chain
complex of $\Tc$ relative to its boundary $\partial \Tc$.

\begin{example}
  We consider the following subdivision $\Tc$ of a rectangle $\Omega$:
  \begin{center}
    \includegraphics{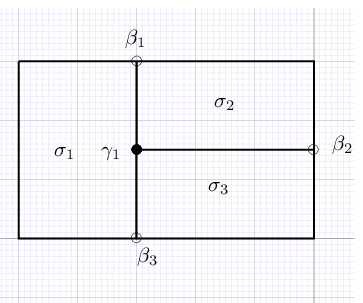}
  \end{center}
  We have
  \begin{itemizedot}
    \item $\partial_2 ([\sigma_1]) = [\gamma_1 \beta_1] + [\beta_3 \gamma_1],
    \partial_2 ([\sigma_2]) = [\gamma_1 \beta_2] + [\beta_1 \gamma_1],
    \partial_2 ([\sigma_3]) = [\gamma_1 \beta_3] + [\beta_2 \gamma_1]$,
    
    \item $\partial_1 ([\beta_1 \gamma_1]) = [\gamma_1], \partial_1 ([\beta_2
    \gamma_1]) = [\gamma_1], \partial_1 ([\beta_3 \gamma_1] = [\gamma_1]$,
    
    \item $\partial_0 ([\gamma_1]) = 0$.
  \end{itemizedot}
  This defines the following complex:
  \[ \begin{array}{lllllllll}
       \mathfrak{R}_{m, m'} (\Tc) : &  & \bigoplus_{i = 1}^3
       [\sigma_i] R_{m, m'} & \rightarrow & \bigoplus_{i = 1}^3 [\beta_i
       \gamma_1] R_{m, m'} & \rightarrow & [\gamma_1] R_{m, m'} & \rightarrow
       & 0
     \end{array} \]
  The matrices of these maps in the canonical (monomial) bases are
  \[ [\partial_2] = \left(\begin{array}{ccc}
       - I & I & 0\\
       0 & - I & I\\
       I & 0 & - I
     \end{array}\right), [\partial_1] = \left(\begin{array}{ccc}
       I & I & I
     \end{array}\right) \]
  where $I$ is the identity matrix of size  $(m + 1) \times (m' +
  1)$ (ie. the dimension of $R_{m, m'}$).

Let us consider the case where $\gamma_{1}= (0,0)$, $(m,m')= (2,2)$
and $\rb$ is the constant distribution $(1,1)$
on $\Tc$. The matrices of the complex $\mathfrak{S}_{2,2}^{1,1} (\Tc)$ are
  \[ [ \bar{\partial}_2] = \left(\begin{array}{ccc}
       - [\Pi_1] & [\Pi_1] & 0\\
       0 & - [\Pi_2] & [\Pi_2]\;\\
       \;[\Pi_3] & 0 & - [\Pi_3]
     \end{array}\right), [\bar{\partial}_1] = \left(\begin{array}{ccc}
       [P_1] & [P_2] & [P_3]
     \end{array}\right) \]
  where $[\Pi_i]$ (resp. $[P_i]$) are the matrices of the projections
$$
\begin{array}{rclrcl}
\Pi_{1}=\Pi_{3}: R_{2,2}&\rightarrow & R_{2,2} /(s^{2}) 
&
\Pi_{2}: R_{2,2}&\rightarrow & R_{2,2} /(t^{2}) \\
p& \mapsto & p \mod s^{2}
&
p& \mapsto & p \mod t^{2}\\
\\
P_{1}=P_{3} : R_{2,2}/ (s^{2})&\rightarrow & R_{2,2} /(s^{2}, t^{2}) 
&
P_{2} : R_{2,2}/(t^{2})&\rightarrow & R_{2,2} /(s^{2}, t^{2}) \\
p \mod s^{2}& \mapsto & p \mod (s^{2},t^{2})
&
p\mod t^{2}& \mapsto & p \mod (s^{2},t^{2}).
\end{array}
$$
The matrices $\Pi_{i}$ are of size $12\times 16$ and the matrices
$P_{i}$ are of size $9\times 12$.
\end{example}

\subsection{Their homology}

In this section, we analyse the homology of the different complexes. The
homology of a chain complex of a triangulation of a (planar) domain is well-known
{\cite{Spanier66}[Chap. 4]},
{\cite{Hatcher02}[Chap. 2]}. 
Since, it is not explicit in the litterature, 
we give in the appendix a simple proof of the exactness of  $\mathfrak{R}_{m, m'}
(\Tc^o)$.

\subsubsection{The 0-homology}
We start by analysing the homology on the vertices.

\begin{lemma}
  $H_0 (\mathfrak{R}_{m, m'} (\Tc^o)) = H_0 (\mathfrak{S}_{m,
  m'}^{\rb} (\Tc^o)) = 0.$
\end{lemma}
\begin{proof} By Proposition \ref{proph0r} in the appendix, we have $H_0 (\mathfrak{R}_{m,
    m'} (\Tc^o)) = 0$. Taking the quotient by $\mathfrak{I}^{\rb}_{m,
  m'} (\tau_i)$, we still get that $\bar{\partial}_1$ is surjective so that 
$H_0 (\mathfrak{S}_{m, m'}^{\rb} (\Tc^o)) = 0$.
\end{proof}

Let us describe in more details $$
H_0 (\mathfrak{I}_{m, m'}^{\rb}
(\Tc^o)) = \bigoplus_{\gamma \in \Tc_0^o}
[\gamma]\mathfrak{I}^{\rb} (\gamma) / \partial_1 ( \bigoplus_{\tau \in
\Tc_1^o} [\tau]\mathfrak{I}^{\rb} (\tau)).
$$ 
We consider the free $R$-module generated by the {\em half-edge}
elements $[\gamma| \tau]$, for all interior edges 
$\tau \in \Tc_1^o$ and all vertices $\gamma \in \tau$.
By convention $[\gamma| \tau] \equiv 0$ if $\gamma \in \partial \Omega$.

For $\gamma \in \Tc_0^o$, let $E_{h} (\gamma)$ (resp. $E_{v} (\gamma)$)
be the set of horizontal (resp. vertical) interior edges
that contain $\gamma$ and let $E (\gamma) = E_{h} (\gamma) \cup E_{v}
(\gamma)$. 
We consider first the map
\begin{eqnarray*}
  \varphi_{\gamma} : \bigoplus_{\tau \in E (\gamma)} [\gamma| \tau]\, R_{(m, m') - \delta (\tau)} & \rightarrow & [\gamma]\,\mathfrak{I}^{\rb}_{m, m'}
  (\gamma)\\\;
  [\gamma| \tau]  & \mapsto & [\gamma]\, \Delta^{\rb} (\tau)
\end{eqnarray*}
By definition of $\mathfrak{I}^{\rb}_{m, m'}(\gamma)$, this map is surjective.
Its kernel is denoted $\mathfrak{K}^{\rb}_{m, m'} (\gamma)$. Let $P_h
(\gamma)$ (resp. $P_v (\gamma)$) be the set of pairs $(\tau,\tau')$ of
distinct horizontal (resp.
vertical) interior edges which contain $\gamma$ (with  $(\tau,\tau')$
identify to $(\tau',\tau)$).
We denote by $P (\gamma) = P_{h} (\gamma) \cup P_{v} (\gamma)$.
 If $\gamma$ is a T-junction,
one of the two sets is empty and the other is a singleton containing one pair.
If $\gamma$ is a crossing vertex, each set is a singleton.

The following proposition gives an explicit description of the kernel
$\mathfrak{K}^{\rb}_{m, m'} (\gamma)$, that we will exploit hereafter.
\begin{proposition} \label{prop:2.3}
  \begin{eqnarray*}
    \mathfrak{K}^{\rb}_{m, m'} (\gamma) & = & \sum_{(\tau, \tau') \in P
    (\gamma)} ([\gamma| \tau] - [\gamma| \tau']) R_{(m, m') - \delta
    (\tau)} \\
   & + & \sum_{\tau \in E_h (\gamma), \tau' \in \text{$E_v (\gamma)$}}
    ([\gamma| \tau]\, \Delta^{\rb} (\tau') - 
    [\gamma| \tau']\, \Delta^{\rb} (\tau)) R_{(m - r - 1, m' - r' - 1)} 
  \end{eqnarray*}
\end{proposition}
\begin{proof}
  Let us suppose first that $\gamma$ is a crossing vertex. We denote by
  $\tau_1, \tau_2$ the horizontal edges, $\tau_3, \tau_4$ the vertical edges
  at $\gamma$. The matrix of the map $\varphi_{\gamma}$ in the basis
  $[\gamma| \tau_i]$ is
  \[ [\varphi_{\gamma}] = \left(\begin{array}{cccc}
       \Delta & \Delta & \Delta' & \Delta'
     \end{array}\right) \]
  where $\Delta = \Delta^{\rb} (\tau_1) = \Delta^{\rb} (\tau_{_2})$, \
  $\Delta' = \Delta^{\rb} (\tau_3) = \Delta^{\rb} (\tau_4) $. Since
  $\Delta$ and $\Delta'$ have no common factor, the kernel of the matrix
  $[\varphi_{\gamma}]$ is generated by the elements $[\gamma| \tau_1] -
  [\gamma| \tau_2]$, $[\gamma| \tau_2]\, \Delta' - [\gamma|
    \tau_3]\, \Delta $, $[\gamma| \tau_3] - [\gamma| \tau_4]$, which give the
  description of $\mathfrak{K}^{\rb}_{m, m'} (\gamma)$ in bi-degree $\leq
  (m,m')$. A similar proof applies when there is no horizontal or
  vertical pair of distinct edges at $\gamma$. This proves the result.
\end{proof}

We use the maps $(\varphi_{\gamma})_{\gamma \in \Tc_0^o}$ to define
\begin{eqnarray*}
\varphi: \bigoplus_{\gamma \in \Tc_0^o} \bigoplus_{\tau
\in E_{\gamma}} [\gamma| \tau]\, R_{(m, m') - \delta (\tau)} & \rightarrow &
\bigoplus_{\gamma \in \Tc_0^o} [\gamma]\,\mathfrak{I}^{\rb}_{m, m'}
(\gamma),
\end{eqnarray*}
so that we have the following exact sequence:
\[ 0 \rightarrow \bigoplus_{\gamma \in \Tc_0^o} \mathfrak{K}^{\rb}_{m, m'} (\gamma) \rightarrow \bigoplus_{\gamma \in \Tc_0^o}
   \bigoplus_{\tau \in E (\gamma)} [\gamma| \tau]\, R_{(m, m') - \delta
   (\tau)} \rightarrow \bigoplus_{\gamma \in \Tc_0^o}
   [\gamma]\,\mathfrak{I}^{\rb}_{m, m'} (\gamma) \rightarrow 0 \]
Using this exact sequence, we can now identify $\bigoplus_{\gamma \in
\Tc_1^o} [\gamma]\,\mathfrak{I}^{\rb}_{m, m'} (\gamma)$ with the
quotient
\[ \bigoplus_{\gamma \in \Tc_1^o} \bigoplus_{\tau \in E_{\gamma}}
   [\gamma| \tau]\, R_{(m, m') - \delta (\tau)} / \sum_{\gamma \in
   \Tc_1^o} \mathfrak{K}^{\rb}_{m, m'} (\gamma) . \]

The next proposition uses this identification  and 
Proposition \ref{prop:2.3} to describe more explicitly 
$H_0 (\mathfrak{I}_{m, m'}^{\rb} (\Tc^o))$:
\begin{proposition}
  We have
  \begin{eqnarray*}
    H_0 (\mathfrak{I}_{m, m'}^{\rb} (\Tc^o)) & = & \bigoplus_{\gamma
    \in \Tc_0^o} \bigoplus_{\tau \in E (\gamma)} [\gamma| \tau]\,
    R_{(m, m') - \delta (\tau)}\\
    & \Bigg \slash & \left( \sum_{(\tau, \tau') \in P (\gamma)} ([\gamma| \tau] -
    [\gamma| \tau']) R_{(m, m') - \delta (\tau)} \right.\\
   && + \sum_{\tau = (\gamma, \gamma') \in \Tc_1^o} ([\gamma|
    \tau] - [\gamma'| \tau]) R_{(m, m') - \delta (\tau)}\\
    && \left. + \sum_{\tau \in E_{h} (\gamma), \tau' \in E_{v} (\gamma)}
    ([\gamma| \tau] \Delta^{\rb} (\tau') - 
    [\gamma| \tau'] \Delta^{\rb} (\tau)) R_{(m, m')  - \delta (\gamma)}  \right) .
  \end{eqnarray*}
\end{proposition}

\begin{proof}
  The application
\begin{eqnarray*}
\partial_1 : \bigoplus_{\tau \in \Tc_1^o} [\tau]\mathfrak{I}^{\rb}_{m, m'} (\tau) &\rightarrow&\bigoplus_{\gamma \in
  \Tc_0^o} [\gamma]\mathfrak{I}^{\rb}_{m, m'} (\gamma)
\end{eqnarray*} 
lift to an application:
\begin{eqnarray*}
    \tilde{\partial}_1 : \bigoplus_{\tau \in \Tc_1^o} [\tau] R_{(m,
    m') - \delta (\tau)} & \rightarrow & \bigoplus_{\gamma \in
    \Tc_0^o}  \bigoplus_{\tau \in E ({\gamma})} [\gamma| \tau] R_{(m,
    m') - \delta (\tau)} \\
    \tau & \mapsto & [\gamma| \tau] - [\gamma'| \tau]
  \end{eqnarray*}
so that the image of $\partial_1$ lift in $\bigoplus_{\gamma \in
  \Tc_0^o}  \bigoplus_{\tau \in E ({\gamma})} [\gamma| \tau] R_{(m,
  m') - \delta (\tau)} $ to
  
  \begin{eqnarray*}
    \tmop{im} \tilde{\partial}_1 & = & \sum_{\tau \in \Tc_1^o}
    ([\gamma| \tau] - [\gamma'| \tau]) R_{(m, m') - \delta (\tau)} .
  \end{eqnarray*}
Consequently,
  \begin{eqnarray*}
    H_0 (\mathfrak{S}_{m, m'}^{\rb} (\Tc^o)) & = & \bigoplus_{\gamma
    \in \Tc_0^o} \bigoplus_{\tau \in E (\gamma)} [\gamma| \tau]\,
    R_{(m, m') - \delta (\tau)} \Bigg \slash \left( \tmop{im} \tilde{\partial}_1 +
    \sum_{\gamma \in \Tc_0^o} \mathfrak{K}^{\rb}_{m, m'} (\gamma)
    \right),
  \end{eqnarray*}
which yields the desired description of $H_0 (\mathfrak{S}_{m, m'}^{\rb} (\Tc))$.
\end{proof}

In the next proposition, we simplify further the description of $H_0 (\mathfrak{I}_{m, m'}^{\rb} (\Tc^o))$:
\begin{proposition}\label{prop:3.5}
$$ 
\begin{array}{lll}
    H_0 (\mathfrak{I}_{m, m'}^{\rb} (\Tc^o)) & = & \oplus_{\rho
    \in \MIS ( {\Tc)}}  \, [\rho]\, R_{(m, m') - \delta
    (\rho)}\\
    & \big \slash & \left( \sum_{\gamma \in \Tc_0^o} (    
[\rho_v(\gamma)]\, \Delta^{\rb}_h (\gamma) - [\rho_h (\gamma)]\, \Delta^{\rb}_v (\gamma) ) R_{(m , m')-\delta (\gamma)}  \right) .
  \end{array}
$$
\end{proposition}

\begin{proof}
  Let $B = \bigoplus_{\gamma \in \Tc_0^o} \bigoplus_{\tau \in
  E ({\gamma})} [\gamma| \tau]\, R_{(m, m') - \delta (\tau)}$, $K = \tmop{im}
  \tilde{\partial}_1 + \sum_{\gamma \in \Tc_0^o} \mathfrak{K}^{r,
  r'}_{m, m'} (\gamma)$ and $$ 
\begin{array}{lll}
    K' & = & \left( \sum_{(\tau, \tau') \in P (\gamma)} ([\gamma| \tau] -
    [\gamma| \tau']) R_{(m, m') - \delta (\tau)}\right.\\
    &  & \left. + \sum_{\tau = (\gamma, \gamma') \in \Tc_1^o}
    ([\gamma| \tau] - [\gamma'| \tau]) R_{(m, m') - \delta (\tau)} \right).
  \end{array}$$
  
  As $K' \subset K \subset B$, we have $B / K \equiv (B / K') / (K / K')$.
  Taking the quotient by $K'$ means, that we identify the horizontal (resp.
  vertical) edges which share a vertex. Thus all horizontal (resp. vertical)
  edges which are contained in a maximal segment $\rho$ of $\Tc$ are
  identify to a single element, that we denote $[\rho]$. As $[\gamma| \tau] =
  0$ if $\gamma \in \partial \Omega$, we also have $[\rho] = 0$ if the maximal segment
  $\rho$ intersects the boundary $\partial\Omega$. This yields the
  desired description of $H_0(\mathfrak{I}_{m, m'}^{\rb} (\Tc))$.
\end{proof}

\begin{definition}
  Let $h_{m, m'}^{\rb} (\Tc) = \dim H_0 (\mathfrak{I}_{m, m'}^{r,
  r'} (\Tc^o)) .$
\end{definition}

\subsubsection{The 1-homology}

We consider now the homology on the edges.
We use the property that  $H_1 (\mathfrak{R}_{m, m'} (\Tc^o))
= 0$ (see Proposition \ref{proph1r} in the appendix).

\begin{proposition}
  $H_1 (\mathfrak{S}_{m, m'}^{\rb} (\Tc^o)) = H_0 (\mathfrak{I}_{m,
  m'}^{\rb} (\Tc^o))$.
\end{proposition}

\begin{proof}
  As $H_0 (\mathfrak{R}_{m, m'} (\Tc^o)) = 0$ and $H_1
  (\mathfrak{R}_{m, m'} (\Tc^o)) = 0$, we deduce from the long exact
  sequence (see Appendix \ref{sec:homology})
  \[ \cdots \rightarrow H_1 (\mathfrak{R}_{m, m'} (\Tc^o)) \rightarrow
    H_1 (\mathfrak{S}_{m, m'}^{\rb} (\Tc^o)) \rightarrow H_0
    (\mathfrak{I}_{m, m'}^{\rb} (\Tc^o)) \rightarrow H_0
    (\mathfrak{R}_{m, m'} (\Tc^o)) \rightarrow \cdots \]
 that $H_1 (\mathfrak{S}_{m, m'}^{\rb} (\Tc^o)) \sim H_0
  (\mathfrak{I}_{m, m'}^{\rb} (\Tc^o))$.
\end{proof}

\subsubsection{The 2-homology}
Finally, the homology on the 2-faces will give us information on the spline
space $\mathcal{S}_{m, m'}^{\rb} (\Tc)$.

We have the following result (proved in Proposition \ref{proph2r} in
the appendix):
\begin{proposition}
  $H_2 (\mathfrak{R}_{m, m'} (\Tc^o)) = R_{m, m'}.$
\end{proposition}
The following proposition relates the spline space $\mathcal{S}_{m, m'}^{\rb} (\Tc)$ 
to an homology module.
\begin{proposition}
  $H_2 (\mathfrak{S}_{m, m'}^{\rb} (\Tc^o)) = \ker \partial_2
  =\mathcal{S}_{m, m'}^{\rb} (\Tc).$
\end{proposition}

\begin{proof}
  An element of $H_2 (\mathfrak{S}_{m, m'}^{\rb} (\Tc^o)) = \ker \overline{\partial}_2$
  is a collection of polynomials $(p_{\sigma})_{\sigma \in
    \Tc_2}$ with $p_{\sigma} \in R_{m, m'}$ and $p_{\sigma} \equiv p_{\sigma'}
  \tmop{mod} \mathfrak{I}_{\tau}^{\rb} (\tau)$ if $\sigma$ and $\sigma'$
  share the (internal) edge $\tau$. By Lemma \ref{lem:reg}, this implies that
  the piecewise polynomial function which is $p_{\sigma}$ on $\sigma$ and
  $p_{\sigma'}$ on $\sigma'$ is of class $C^{\rb}$ across the edge $\tau$.
  As this is true for all interior edges, $(p_{\sigma})_{\sigma \in  \Tc_2} \in \ker \underline{\partial}_2 $ is a piecewise polynomial function of
  $R_{m, m'}$ which is of class $C^{\rb}$, that is an element of
  $\mathcal{S}_{m, m'}^{\rb} (\Tc)$.
\end{proof}


 
\section{Lower and upper bounds on the dimension}\label{sec:3}
In this section, are the main results on the dimension of the spline
space $\mathcal{S}_{m, m'}^{\rb} (\Tc)$.
\begin{theorem}
 \label{thm:dim}
Let $\Tc$ be a T-mesh and let
 $\rb$ be a smoothness distribution on $\Tc$. Then
\begin{eqnarray}\label{eq:3.1}
{\dim \mathcal{S}_{m, m'}^{\rb} (\Tc)} 
& = & \sum_{\sigma\in \Tc_{2}\ } (m + 1) (m' + 1) \label{eq:dim}  \\
& -  & \sum_{\tau\in \Tc_{1}^{o,h}} (m + 1) (\rb (\tau) + 1)  -
\sum_{\tau\in \Tc_{1}^{o,v}} (m' + 1) (\rb (\tau) + 1) \nonumber\\
&+& \sum_{\gamma\in \Tc_{0}^{o}}(\rb_{h} (\gamma)  + 1) ( \rb_{v} (\gamma) + 1) \nonumber\\
   & + & h_{m, m'}^{\rb} (\Tc). \nonumber
 \end{eqnarray}
where $h_{m, m'}^{\rb} (\Tc) = \dim H_0 (\mathfrak{I}_{m,
  m'}^{\rb} (\Tc^o))$.
\end{theorem}
\begin{proof}
 The complex
{\small  \[ \mathfrak{S}_{m, m'}^{\rb} (\Tc^o) : \oplus_{\sigma \in
    \Tc_2} [\sigma] R_{m, m'} \longrightarrow \oplus_{\tau \in
    \Tc_1^o} [\tau] R_{m, m'} /\mathfrak{I}^{\rb}_{m, m'} (\tau)
    \longrightarrow \oplus_{\gamma \in \Tc_0^o} [\gamma] R_{m, m'}
    /\mathfrak{I}^{\rb}_{m, m'} (\gamma) \longrightarrow 0 \]
}
 induces the following relations
 \[ \dim (\oplus_{\sigma \in \Tc_2} [\sigma] R_{m, m'}) - \dim
    (\oplus_{\tau \in \Tc_1^o} [\tau] R_{m, m'} /\mathfrak{I}^{r,
    r'}_{m, m'} (\tau)) + \dim (\oplus_{\gamma \in \Tc_0^o} [\gamma]
    R_{m, m'} /\mathfrak{I}^{\rb}_{m, m'} (\gamma)) \]
 \[ = \dim (H_2 (\mathfrak{S}_{m, m'}^{\rb} (\Tc^o))) - \dim (H_1
    (\mathfrak{S}_{m, m'}^{r, r} (\Tc^o))) + \dim (H_0
    (\mathfrak{S}_{m, m'}^{\rb} (\Tc^o))) \]
 As $H_2 (\mathfrak{S}_{m, m'}^{\rb} (\Tc^o)) =\mathcal{S}_{m,
 m'}^{\rb} (\Tc)$, $H_0 (\mathfrak{S}_{m, m'}^{\rb}
 (\Tc^o)) = 0$ and $H_1 (\mathfrak{S}_{m, m'}^{\rb}
 (\Tc^o)) = H_0 (\mathfrak{I}^{\rb}_{m, m'} (\Tc^o))$, we
 deduce that
$$ 
\begin{array}{lll}
   \dim \mathcal{S}_{m, m'}^{\rb} (\Tc) & = & \dim (\oplus_{\sigma
   \in \Tc_2} [\sigma] R_{m, m'}) - \dim (\oplus_{\tau \in
   \Tc_1^o} [\tau] R_{m, m'} /\mathfrak{I}^{\rb}_{m, m'} (\tau))\\
   & + & \dim (\oplus_{\gamma \in \Tc_0^o} [\gamma] R_{m, m'}
   /\mathfrak{I}^{\rb}_{m, m'} (\gamma)) + \dim (H_0 (\mathfrak{I}^{r,
   r'}_{m, m'} (\Tc)))
 \end{array}
$$
which yields the dimension formula \eqref{eq:3.1} using Lemma \ref{lem:dim}.
\end{proof}

As an immediate corollary of this theorem and of Proposition
\ref{prop:3.5}, we deduce the following result: 
\begin{corollary}
 \label{cor:dim:exact} If the T-mesh $\Tc$ has no maximal
interior segments then $h_{m, m'}^{\rb} (\Tc)=0$.
\end{corollary}

In the case of a constant smoothness distribution, Theorem \ref{thm:dim}
is written as follows:
\begin{theorem} \label{thm:dim:cst}
Let $\Tc$ be a T-mesh and let
 $\rb=(r,r')$ be a constant smoothness distribution on $\Tc$. Then
\begin{eqnarray}
{\dim \mathcal{S}_{m, m'}^{\rb} (\Tc)} 
& = & (m + 1) (m' + 1) f_2 \label{eq:dim:cst}  \\
& -  & \left((m + 1) (r' + 1) f_1^h + (m' + 1) (r + 1) f_1^v \right) \nonumber\\
&+& (r + 1) (r' + 1) f_0 \nonumber\\
   & + & h_{m, m'}^{\rb} (\Tc). \nonumber
 \end{eqnarray}
 where
  \begin{itemizedot}
   \item $f_2$ is the number of 2-faces $\in \Tc_2$,
    
    \item $f_1^h$ (resp. $f_1^v$) is the number of horizontal (resp. vertical)
    interior edges $\in \Tc_1^o$,
    
    \item $f_0$ is the number of interior vertices $\in \Tc_0^o$,
    
    \item $h_{m, m'}^{\rb} (\Tc) = \dim H_0 (\mathfrak{I}_{m,
  m'}^{\rb} (\Tc^o)) .$
 \end{itemizedot}
\end{theorem}

We are going now to bound $h_{m, m'}^{\rb} (\Tc)$ for general T-meshes.

\begin{definition}
Let $\iota$ be an ordering of $\MIS(\Tc)$ that is a map from $\MIS(\Tc)$ to $\mathbb{N}$. 
For $\rho \in \MIS (\Tc)$,
let $\Gamma_{\iota}(\rho)$ be the set of vertices $\gamma$ of $\rho$
which are not on a maximal interior segment $\rho'\in \MIS (\Tc)$ with
$\iota (\rho') > \iota (\rho)$. 
The number of elements of 
$\Gamma_{\iota}(\rho)$
is denoted 
$\lambda_{\iota} (\rho)$.
\end{definition}

We define now the weight of a maximal interior segment.
\begin{definition}
For $\rho\in \MIS (\Tc)$, let 
\begin{itemize}
 \item $\omega_{\iota} (\rho) = \sum_{\gamma\in  \Gamma_{\iota}(\rho)}
  (m-\rb_{v} (\gamma))$ if $\rho\in \MIS_{h} (\Tc)$.
 \item $\omega_{\iota} (\rho) = \sum_{\gamma\in  \Gamma_{\iota}(\rho)} (m'-\rb_{h} (\gamma))$ if $\rho\in \MIS_{v} (\Tc)$.
\end{itemize}
We called it the {\em weight} of $\rho$.
\end{definition}
As in the usual spline terminology, for an interior point $\gamma\in
\Tc_{0}^{o}$, we call 
$\gamma$ $m-\rb_{h} (\gamma)$ (resp. $m'-\rb_{h} (\gamma)$ )
the horizontal (resp. vertical) multiplicity of $\gamma$.

If $\rho$ is horizontal (resp. vertical), the weight of $\rho$ is 
the sum of the vertical (resp. horizontal) multiplicities of the
vertices $\gamma \in \Gamma_{\iota} (\rho)$. 

Notice that if $\rb= (r,r')$ is a constant smoothness distribution on
$\Tc$, then $\omega_{\iota} (\rho)= (m -r) \lambda (\rho)$ for $\rho
\in \MIS_{h} (\Tc)$
(resp. $\omega_{\iota} (\rho)= (m' -r') \lambda (\rho)$ for $\rho \in \MIS_{v} (\Tc)$).
\begin{example} We consider the T-mesh of Example \ref{ex:1.9} with $(m,m') =
  (2,2)$, the constant smoothness distribution $\rb= (1,1)$ and the ordering
of the maximal interior segments $\iota (\rho_{i})=i$ for $i=1, \ldots, 4$. Then we have
\begin{itemize}
 \item $\omega_{\iota} (\rho_{1}) = (2-1) \times 2  =2$,
 \item $\omega_{\iota} (\rho_{2}) = (2-1) \times 2  =2$,
 \item $\omega_{\iota} (\rho_{3}) = (2-1) \times 3  =3$,
 \item $\omega_{\iota} (\rho_{4}) = (2-1) \times 3  =3$,
\end{itemize}
since the multiplicity of a vertex is $2-1=1$ and the interior points of
$\rho_{1}, \rho_{2}$ are
not in $\Gamma_{\iota} (\rho_{1})$ or $\Gamma_{\iota} (\rho_{2})$.
\end{example}
 
In the following, we will drop the index $\iota$ to simplify the 
notations, assuming that the ordering $\iota$ is fixed. 
Then we have the following theorem:

\begin{theorem}\label{thm:h}\label{thm:dim:bound}
Let $\Tc$ be a T-mesh and let $\rb$ be a smoothness distribution on $\Tc$. Then
  \begin{eqnarray*}
   0\ \leq\ h_{m, m'}^{\rb} (\Tc) & \leq & \sum_{\rho \in \MIS_h
    (\Tc)}
(m + 1 - \omega (\rho))_+ \times (m' - \rb(\rho))\\
    & + & \sum_{\rho \in \MIS_v (\Tc)}  (m - \rb (\rho)) \times \left( m' +
   1 - \omega (\rho) )\right)_+.
  \end{eqnarray*}
\end{theorem}
\begin{proof} Let $\rho_{1},\ldots, \rho_{l}$ be the maximal interior
  segments of $\Tc$.
 By Proposition \ref{prop:3.5}, $h_{m, m'}^{\rb} (\Tc)$
 is the dimension of the
 quotient in bi-degree $\leq (m, m')$ of the module $M \assign \oplus_{i= 1}^l [\rho_i]\, R$ by the
 module $K$ generated by the following relations:
 for each vertex $\gamma \in \Tc_0^{o}$ which is on a maximal interior
 segment,
 \begin{itemizedot}
   \item $\Delta^{\rb} (\rho_j) [\rho_{_i}] - \Delta^{\rb} (\rho_i) [\rho_j]$ if $\gamma$ is the
   intersection of the maximal interior segments $\rho_i$ and $\rho_j$,
   
   \item $\Delta^{\rb} (\rho) [\rho_{_i}]$ if $\gamma$ is the intersection of the
   maximal interior segment $\rho_i$ with another maximal segment
   $\rho$ intersecting $\partial \Omega$.
 \end{itemizedot}
To compute the dimension of $M/K$ in bi-degree $\leqslant (m,m')$, 
we use a graduation on $M$ given by the indices of the segments.
For $r:=\sum_i p_i [\rho_i] \in M$ 
(with $p_i \in R_{(m, m') - \delta (\rho_i)}$),
 let $\In (r)$ be the element $p_{i_0} [\rho_{i_0}]$
 where $i_0$ is the maximal index such that $p_i \neq 0$. 
We denote it by $\In (r)$ and called it the initial of $r$. Let $\In (K)=\{ \In (k)\mid k\in K \}$.
The dimension
 $h_{m, m'}^{\rb} (\Tc)$ is then
 \[ h_{m, m'}^{\rb} (\Tc) = \dim (M / K) = \dim (M / \In (K)). \]
Notice that $\In (K)$ contains the multiples in bi-degree $\leqslant (m, m')$ of
 \begin{itemizedot}
   \item $\Delta^{\rb} (\rho_j) [\rho_{_i}]$ if $\gamma$ is the intersection of two
   maximal interior segments $\rho_i$ and $\rho_j$ with $i > j$,
   
   \item $\Delta^{\rb} (\rho) [\rho_{_i}]$ if $\gamma$ is the intersection of the
   maximal interior segment $\rho_i$ with a maximal segment $\rho$
   intersecting $\Omega$.
 \end{itemizedot}
Let $L_i$ be the vector space spanned by these
 initials in bi-degree $\leqslant (m, m')$, which are multiples of
 $[\rho_{i}]$.  
By definition,  for each $\gamma\in \Gamma (\rho_{i})$, we have a
generator $\Delta^{\rb} (\rho)[\rho_{i}]$ in $L_{i}$ for
$\{\gamma\}=\rho_{i}\cap \rho$.
By Proposition \ref{dim:apolar}, $L_i$ is of dimension
\begin{itemizedot}
 \item $\min (m + 1, \omega (\rho_{i})) \times (m' - \rb (\rho_{i}))$ 
if $\rho_i \in  \MIS_h (\Tc)$, 
   
 \item $\min (m' + 1, \omega (\rho_{i})) \times (m - \rb (\rho_{i}))$
 if $\rho_i \in \MIS_v (\Tc)$. \end{itemizedot}
Thus the dimension of $[\rho_{i}] R_{(m, m') - \delta (\rho_i)} / L_i$ is
\begin{itemizedot}
   \item 
$(m + 1 - \omega (\rho_{i}))_+ \times (m' - \rb(\rho_{i}))$
if $\rho_i \in  \MIS_h (\Tc)$,
   
   \item  $(m - \rb (\rho_{i})) \times ( m' +  1 - \omega (\rho_{i}))_+$
 if $\rho_i \in \MIS_v (\Tc)$.
\end{itemizedot}
 As $\tmop{In} (K) \supset \sum_i L_i$, we have
\begin{eqnarray*}
h_{m, m'}^{\rb} (\Tc) &=& \dim (M / \tmop{In} (K))\\ 
&\leq& \dim ( [\rho_{_i}] \, R_{(m, m') -
    \delta (\rho_i)}/  (\sum_i L_i)) = \sum_i \dim \left( [\rho_{_i}] \,R_{(m, m') - \delta (\rho_i)}/ L_i\right). 
\end{eqnarray*}
This gives the announced bound on $h_{m, m'}^{\rb} (\Tc)$,
 using the previous computation of $\dim( [\rho_{_i}]\, R_{(m, m') - \delta (\rho_i)}/ L_i)$.
\end{proof}

\begin{definition}
The T-mesh $\Tc$ with a smoothness distribution $\rb$ is $(k,k')$-weighted if 
\begin{itemize}
\item $\forall \rho \in \MIS_h (\Tc)$, $\omega (\rho)\geq k$ 
\item $\forall \rho \in \MIS_v (\Tc)$, $\omega (\rho) \geq k'$; 
\end{itemize}
\end{definition}
 
\begin{theorem}\label{thm:weighted}
Let $\Tc$ be a T-mesh and let $\rb$ be a smoothness distribution on
$\Tc$.  If $\Tc$ is $(m+1,m'+1)$-weighted, then $h_{m, m'}^{\rb} (\Tc) = 0$.
\end{theorem}
\begin{proof} 
By definition, $\forall \rho \in \MIS_h (\Tc)$, $ \omega
(\rho) \geq m + 1$ 
(i.e. $(m+1 - \omega (\rho))_{+}=0$)
and
$\forall \rho \in \MIS_v (\Tc)$, 
$\omega (\rho) \geq m' + 1$ 
(i.e. $(m'+1 - \omega (\rho))_{+}=0$). 
By Theorem \ref{thm:h}, we directly deduce that $h_{m, m'}^{\rb} (\Tc) = 0$.
\end{proof}
 
Here is a direct corollary which generalizes a result in \cite{LiWaZha06}:
\begin{corollary}
Suppose that the end points of a maximal interior segment $\rho\in \MIS (\Tc)$ are in
$\Gamma (\rho)$.
If for each horizontal (resp. vertical) maximal interior segment
$\rho\in \MIS (\Tc)$, the
sum of the vertical (resp. horizontal) multiplicities of the end
points and of the vertices of $\rho$ on a maximal segment connected to the
boundary is greater than or equal to $m+1$ (resp. $m'+1$), then $h_{m, m'}^{\rb} (\Tc)=0$.
\end{corollary}
\begin{proof}
Let $\rho$ be a maximal interior segment of $\Tc$.
By hypothesis, the end points of $\rho$ are in $\Gamma (\rho)$.
As any point $\gamma \in \rho$ which is also on a maximal segment connected to the
boundary is in $\Gamma (\rho)$, the hypothesis implies that $\omega (\rho)\ge m+1$ if $\rho$ is
horizontal and $\omega (\rho)\ge m'+1$ if $\rho$ is vertical.
We deduce by Theorem \ref{thm:dim:bound} that $h_{m, m'}^{\rb} (\Tc)=0$.
\end{proof}

Another case where $h_{m, m'}^{\rb} (\Tc)$ is known is described in the
next proposition.
\begin{proposition}\label{prop:upper}
If $\forall \rho \in \MIS_h (\Tc)$, $\omega (\rho) \leq m + 1$ and $\forall
   \rho \in \MIS_v (\Tc)$, $\omega (\rho) \leq m' + 1$,
   then
\begin{eqnarray}\label{eq:upper}
    h_{m, m'}^{\rb} (\Tc) & = & \sum_{\rho \in \MIS_h
    (\Tc)}
(m + 1 - \omega (\rho))_+ \times (m' - \rb(\rho))\\
    & + & \sum_{\rho \in \MIS_v (\Tc)}  (m - \rb (\rho)) \times \left( m' +
   1 - \omega (\rho) )\right)_+.\nonumber
\end{eqnarray}

\end{proposition}
\begin{proof}
In the case where 
$\forall \rho \in \MIS_h (\Tc)$, $\omega (\rho) \leq m + 1$ and $\forall
   \rho \in \MIS_v (\Tc)$, 
$\omega (\rho) \leq m' + 1$, 
Proposition \ref{dim:apolar} implies that there is no relations in bi-degree $\leqslant (m, m')$ of the
monomial multiples of $\Delta^{\rb} (\rho_j) [\rho_{_i}]$, $\Delta^{\rb} (\rho) [\rho_{_i}]$ for
$i = 1, \ldots, l, j<i$ using the same notation as in the proof of Theorem \ref{thm:dim:bound}. 
This implies that $\In (K) = \oplus_i L_i$, which
shows that 
$$ 
h_{m, m'}^{\rb} (\Tc) = \dim (M / \In (K)) = 
\sum_i \dim (R_{(m, m') - \delta (\rho_i)}
    [\rho_{_i}] / L_i).
$$
Thus the equality \eqref{eq:upper} holds.
\end{proof}

As a corollary, we have the following result for constant smoothness distribution:
\begin{theorem}\label{thm:h:cst} 
  \label{thm:dim:bound:cst}
Let $\Tc$ be a T-mesh and let
  $\rb= (r,r')$ be a constant smoothness distribution on $\Tc$. Then
  \begin{eqnarray*}
0 \leq    h_{m, m'}^{\rb} (\Tc) & \leq & \sum_{\rho \in \MIS_h
    (\Tc)} (m + 1 - (m - r) \lambda (\rho))_+ \times (m' - r')\\
    & + & \sum_{\rho \in \MIS_v (\Tc)} (m - r) \times \left( m' +
    1 - (m' - r') \lambda (\rho) \right)_+ .
  \end{eqnarray*}
Moreover, equality holds in the following cases:
  \begin{itemizedot}
    \item $\forall \rho \in \MIS_h (\Tc)$, $(m - r) \lambda (\rho)
    \geq m + 1$ and $\forall \rho \in \MIS_v (\Tc)$, $(m' -
    r') \lambda (\rho) \geq m' + 1$; 
    
    \item $\forall \rho \in \MIS_h (\Tc), (m - r) \lambda (\rho)
    \leq m + 1$ and $\forall \rho \in \MIS_v (\Tc)$, $(m' -
    r') \lambda (\rho) \leq m' + 1$.
  \end{itemizedot}
\end{theorem}

\section{Hierarchical T-meshes} \label{sec:4}

We consider now a special family of T-meshes, which can be defined
by recursive subdivision from an initial rectangular domain $\Omega$. Their study is motivated
by practical applications, where local refinement of tensor-product spline spaces are
considered eg. in isogeometric analysis \cite{hughes:CMAME2005}.
\begin{definition}
 A hierarchical T-mesh is either the initial axis-aligned rectangle $\Omega$ or obtained from a
 hierarchical T-mesh by splitting a cell along a vertical or
 horizontal line.
\end{definition}

A hierarchical T-mesh will also be called a T-subdivision. 
It can be represented by a subdivision tree where the
nodes are the cells obtained during the subdivision and the children of a cell
$\sigma$ are the cells obtained by subdividing $\sigma$. 

In a hierarchical T-mesh $\Tc$, the maximal interior segments are
naturally ordered in the way they appear during the subdivision process. 
This is the ordering $\iota$ that we will consider hereafter. 
Notice that a maximal interior segment $\rho$ is transformed by
a split either into a maximal segment which intersects
$\partial\Omega$ or into a larger maximal interior segment with a
larger weight. 

We remark that in a hierarchical T-mesh, if $\rho_{i}$ is blocking
$\rho_{j}$ then $i< j$.

\begin{example}
 Here are an hierarchical T-mesh (case a) and a non-hierarchical
 T-mesh (case b):  
 \begin{center}
   \includegraphics[height=5cm]{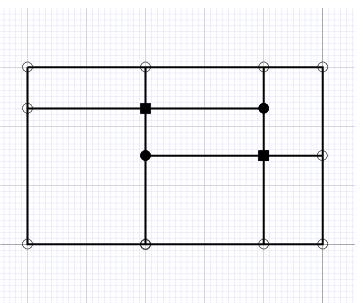}
   \includegraphics[height=5cm]{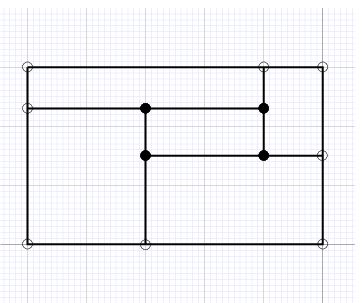}\\
(a) \ \hspace{3cm}\ \ \ \ \ \ \ \ \  (b)
 \end{center}
 
  
\end{example}

\subsection{Dimension formula for hierarchical T-meshes}

As a corollary of Theorem \ref{thm:dim:bound}, we deduce the following result, also proved in
{\cite{deng06}, \cite{LiWaZha06}, \cite{HuDeFeChe06}}:

\begin{proposition}
\label{prop:dim:exact} Let $\Tc$ be a hierarchical T-mesh and
let $\rb= (r,r')$ be a constant smoothness distribution on $\Tc$. For
$m \geq 2 r + 1$ and $m' \geq 2 r' + 1$, we have 
$h_{m, m'}^{r,r'} (\Tc) = 0.$
\end{proposition}

\begin{proof}
We order the maximal interior segments in the way they appear during the
subdivision. If a segment $\rho_i \in \MIS (\Tc)$ is blocking $\rho_j
\in \MIS (\Tc)$, we must have $i < j$. This shows that the end
points of $\rho_i$ are in $\Gamma (\rho_i)$. Thus, $\lambda (\rho_i)
\geq 2$.

As \ $m \geq 2 r + 1$, we have
\[ (m - r) \lambda (\rho_i) \geq 2 (m - r) \geq m + (m - 2 r)
   \geq m + 1. \]
Thus, $(m + 1 - (m - r) \lambda (\rho_i))_+ = 0$. Similarly $(m' + 1 - (m' -
r') \lambda (\rho_i))_+ = 0$ holds since $m' \geq 2 r' + 1$. By Theorem
\ref{thm:h}, we deduce that $h_{m, m'}^{\rb} (\Tc) = 0$.
\end{proof}

Theorem \ref{thm:weighted} leads us to the following construction rule
of a T-subdivision $\Tc$ for which $h_{m, m'}^{\rb} (\Tc) = 0.$
\begin{algorithm}[$(k,k')$-weighted subdivision rule]\ \\
For each 2-face $\sigma$ of a T-mesh to be subdivided,
\begin{enumerate}
 \item Split $\sigma$ with the new edge $\tau$;
 \item If the edge $\tau$ does not extend an existing segment, extend $\tau$
(on one side and/or the other) so that the maximal segment containing 
$\tau$ is either intersecting $\partial \Omega$ or horizontal
(resp. vertical) and of weight $\geq k$ (resp. $\geq k'$).
\end{enumerate}
\end{algorithm}
If such a rule is applied in the construction of a T-subdivision,
\begin{itemize}
 \item either a new maximal interior segment is constructed so that
its weight is $\ge k$ if it is a horizontal maximal interior segment 
(resp. $\ge k'$ it it is a vertical maximal interior segment),
\item or an existing maximal interior segment is extended
  and its weight is also increased,
\item or a maximal segment intersecting $\partial \Omega$ is constructed.
\end{itemize}
In all cases, if we start with a $(k,k')$-weighted T-mesh, we obtain a
new T-mesh, which is also $(k,k')$-weighted.

By Theorem \ref{thm:weighted}, if $k \ge m+1$ and $k'\ge m'+1$ then
$h_{m, m'}^{\rb} (\Tc) = 0$ and the dimension of 
$\dim \mathcal{S}_{m, m'}^{\rb} (\Tc)$, given by formula \eqref{eq:dim}, depends
  only on the number of cells, interior segments and interior vertices of
  $\Tc$. From this analysis, we deduce the dimension formula of the
  space of Locally Refined splines described in \cite{DoLyPe10}.
\section{Examples}
In this section, we analyse the dimension formula of spline spaces of small bi-degree and
small constant smoothness distribution $\rb= (r,r')$ on a T-mesh $\Tc$.

\subsection{Bilinear $C^{0,0}$ T-splines}

We consider first piecewise bilinear polynomials on $\Tc$, which are
continuous, that is $m = m' = 1$ and $r = r' = 0$. By Proposition
\ref{prop:dim:exact}, we have $h_{1, 1}^{0, 0} (\Tc) = 0$. Using
Theorem \ref{thm:dim} and Lemma \ref{lem:nbf} in the appendix, we obtain:
\begin{equation}
  \dim \mathcal{S}_{1, 1}^{0, 0} (\Tc) = 4 f_2 - 2 f_1^o + f_0^o =
  f_0^+ + f_0^b. \label{eq:dim:c0}
\end{equation}

\subsection{Biquadratic $C^{1,1}$ T-splines}
Let us consider now the set of piecewise biquadratic functions on a
T-mesh $\Tc$, which are $C^{1}$.
For $m = m' = 2$ and $r = r' = 1$, Theorem \ref{thm:dim} and Lemma
\ref{lem:nbf} again yield 
\begin{equation}
  \dim \mathcal{S}_{2, 2}^{1, 1} (\Tc) = 9 f_2 - 6 f_1^o + 4 f_0^o +
  h_{2, 2}^{1, 1} (\Tc) = f_0^+ - \frac{1}{2} f_0^T + \frac{3}{2}
  f_0^b + 3 + h_{2, 2}^{1, 1} (\Tc) .
\end{equation}
If the T-mesh $\Tc$ is $(3,3)$-weighted, then by Theorem \ref{thm:weighted}, we have a  $h_{2, 2}^{1, 1}(\Tc)=0$, but this is not always the case.
\begin{example}\label{ex:5.1}
  Here is an example where $h_{2, 2}^{1, 1} (\Tc) = 1$ by Proposition
  \ref{prop:upper}, since there is one maximal interior segment $\rho$
  with $\omega (\rho) = 2-1 + 2-1 = 2$:
  
  \begin{center}
    \includegraphics{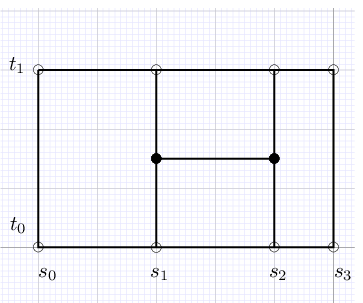}
  \end{center}
  
  The dimension of $\mathcal{S}_{2, 2}^{1, 1} (\Tc)$ is $9 \times 4 -
  6 \times 5 + 4 \times 2 + h_{2, 2}^{1, 1} (\Tc) = 14 + 1 = 15$.
  Notice that the dimension is the same without the (horizontal) interior
  segment. Thus a basis of \ $\mathcal{S}_{2, 2}^{1, 1} (\Tc)$ is the
  tensor product B-spline basis corresponding to the nodes $s_0, s_0, s_0, s_1,
  s_2, s_3, s_3, s_3$ in the horizontal direction and the nodes $t_0, t_0,
  t_0, t_1, t_1, t_1$ in the vertical direction.
\end{example}

\begin{example}\label{ex:5.2}
  Here is another example. We subdivide the T-mesh $\Tc_1$ to
  obtain the second T-mesh $\Tc_2$:
  \begin{center}
    \includegraphics{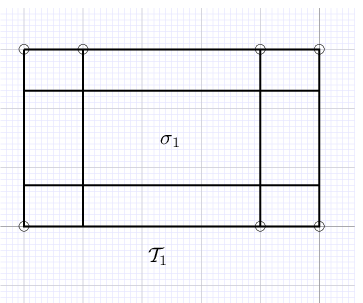} \ 
    \includegraphics{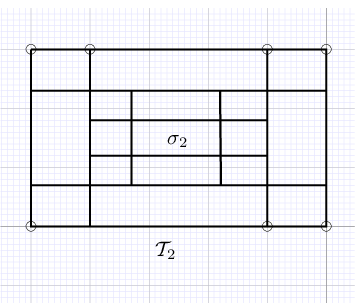}
  \end{center}
  
  Doing this, we increase the number of cells by $9 - 1 = 8$, the number of
  interior edges by $24 - 4 = 20$, the number of interior points by $16 - 4 =
  12$. The dimension of the spline space increases by $9 \times 8 - 6 \times
  20 + 4 \times 12 + h_{2, 2}^{1, 1} (\Tc_2) - h_{2, 2}^{1, 1}
  (\Tc_1)= h_{2, 2}^{1, 1} (\Tc_2) - h_{2, 2}^{1, 1}
  (\Tc_1)$. Since there is no
  maximal interior segment in $\Tc_{1}$, by Corollary \ref{cor:dim:exact} we have $h_{2,
  2}^{1, 1} (\Tc_1) = 0$. Choosing a proper ordering of the interior
  segments, we deduce by Theorem \ref{thm:dim} that $h_{2, 2}^{1, 1}
  (\Tc_2) \leq 1$. Suppose that $\sigma_1 = [a_0, a_3] \times
  [b_0, b_3]$ and $\sigma_2 = [a_1, a_2] \times [b_1, b_2]$. 
For $u_{0}\le u_{1}\le u_{2}\le u_{3}\in \RR$, let $N (u;u_{0},
u_{1},u_{2},u_{3})$ be the B-spline basis function in the variable $u$
of degree $2$ associated to the
nodes $u_{0}, \ldots, u_{3}$ (see \cite{deBoor01}). 
Then the piecewise polynomial function
  \[ N (s ; a_0, a_1, a_2, a_3) \times N (t ; b_0, b_1, b_2, b_3) \]
  is an element of $\mathcal{S}_{2, 2}^{1, 1} (\Tc_2)$, with support
  in $\sigma_1$. It is not in $\mathcal{S}_{2, 2}^{1, 1} (\Tc_1)$,
  since the function is not polynomial on $\sigma_1$. Thus we have $\dim
  \mathcal{S}_{2, 2}^{1, 1} (\Tc_2) = \dim \mathcal{S}_{2, 2}^{1, 1}
  (\Tc_1) + 1$. 

Notice that $\Tc_2$ is $(2,2)$-weighted but not $3$-weighted, since 
any new maximal segment intersects two of the other new maximal segments.
\end{example}
For a general hierarchical T-mesh, we consider a sequence of
T-meshes $\Tc_{0}, \ldots \Tc_{l}$ where $\Tc_{0}$ has one cell, 
$\Tc_{l}=\Tc$ and such that $\Tc_{i+1}$ refines $\Tc_{i}$ by inserting new edges.
We can assume that at each level $i\neq 0$ a new maximal interior segment $\rho_{i}$ appears
and that we number the maximal interior segments of $\Tc$ in
the order they appear during this subdivision. 
Notice that any maximal interior segment of $\Tc$ extends one
of the maximal interior segments $\rho_{i}$ and thus its weight is bigger.
Notice also that the maximal interior segment $\rho_{i}$ introduced at
level $i$ extends to a maximal segment
of $\Tc$, which may intersect the boundary. In this case, it is not involved in
the dimension upper bound. Then, we have the following corollary:
\begin{proposition} Let $\Tc$ be a hierarchical T-mesh.
\begin{equation}\label{eq:5.3}
9 f_2 - 6 f_1^o + 4 f_0^o \le \dim \mathcal{S}_{2, 2}^{1, 1} (\Tc) \le 9 f_2 - 6 f_1^o + 4 f_0^o + \sigma
\end{equation} 
where $\sigma$ is the number of levels of the subdivision where a maximal interior segment  with
no-interior point is introduced.
\end{proposition}
\begin{proof}
Consider the new maximal interior segment $\rho_{i}$ of $\Tc_{i}$ appearing at level $i$.
By construction, the end points of $\rho_{i}$ are not on maximal interior segments of bigger
index. 
Thus $\omega (\rho_{i})\ge 2$. If $\rho_{i}$ contains an interior vertex then $\omega (\rho_{i})\ge 3$
and $(3 - \omega (\rho_{i}))_{+}=0$. Otherwise $(3 - \omega (\rho_{i}))_{+}=1$.
As $\rho_{i}$ extends to a maximal segment
$\tilde{\rho}_{i}$ which is either interior or intersecting the boundary, we
have $(3 - \omega (\tilde{\rho}_{i}))_{+} \le (3- \omega
(\rho_{i}))_{+}$ with the convention
that $\omega (\tilde{\rho}_{i})=3$ if $\tilde{\rho}_{i}$ is
intersecting the boundary.

By Theorem \ref{thm:dim:bound}, we have $$ 
0 \le h_{2,2}^{1,1} (\Tc)
\le \sum_{i=1}^{l} (3 - \omega (\tilde{\rho}_{i}))_{+}
\le \sum_{i=1}^{l} (3 - \omega (\rho_{i}))_{+}
$$ 
where $l$ is the number of levels in the subdivision,
$\rho_{i}$ is the maximal interior segment introduced at level $i$ and 
$\tilde{\rho}_{i}$ is its extension in $\Tc$.
By the previous remarks, $\sum_{i=1}^{l} (3 - \omega
(\rho_{i}))_{+}=\sigma$ is the number of levels of the subdivision where a maximal interior segment  with
no interior point is introduced.
Using Theorem \ref{thm:dim}, this proves the bound on the dimension of $\mathcal{S}_{2, 2}^{1, 1}
(\Tc)$. 
\end{proof}
Examples \ref{ex:5.1} and \ref{ex:5.2} show that the dimension
can be given by the upper bound. 
On other other hand, for any (3,3)-weighted  hierarchical
subdivision, the lower bound is reached. This shows that the
inequalities \eqref{eq:5.3} are optimal for 
$\dim \mathcal{S}_{2, 2}^{1, 1} (\Tc)$.
\begin{remark}
In the case of a hierarchical subdivision where some cells of a given
level are subdivided  (as $\sigma_{1}$ in Example \ref{ex:5.2}) into 9 sub-cells
which have the same length and height, it can be proved that the dimension
is in fact:
$$ 
\dim \mathcal{S}_{2, 2}^{1, 1} (\Tc) = 9 f_2 - 6 f_1^o + 4 f_0^o + \sigma
$$
where $\sigma$ is the number of isolated subdivided cells (i.e. the
cell is subdivided, not touching the boundary and the adjacent cells sharing an
edge are not subdivided) at some level of the subdivision. Indeed, any
maximal interior segment subdividing a non-isolated cell contains an
interior point and is not involved in the upper bound. As in Example
\ref{ex:5.2}, only the isolated cells have a maximal interior segment
with no interior points. This example 
also shows that a new $C^{(1,1)}$ bi-quadratic basis element can be
constructed for each isolated cell, proving that the upper bound is
reached. This gives a dimension formula similar to the one in
\cite{Deng08TSPL22}, for a slightly different subdivision strategy.
\end{remark}

\subsection{Bicubic $C^{1,1}$ T-splines}

For $m = m' = 3$ and $r = r' = 1$, that is for piecewise bicubic polynomial
functions which are $C^{1}$, Proposition \ref{prop:dim:exact} yields $h_{3,
  3}^{1, 1} (\Tc) = 0$. Using Theorem \ref{thm:dim} and Lemma
\ref{lem:nbf} in the appendix, we obtain:
\begin{equation}
  \dim \mathcal{S}_{3, 3}^{2, 2} (\Tc) = 16 f_2 - 8 f_1^o + 4 f_0^o =
  4 (f_0^+ + f_0^b). \label{eq:dim:d3c1}
\end{equation}

\subsection{Bicubic $C^{2,2}$ T-splines}

For $m = m' = 3$ and $r = r' = 2$, by Theorem \ref{thm:dim} and Lemma
\ref{lem:nbf}, we have:
\begin{equation}
  \dim \mathcal{S}_{3, 3}^{2, 2} (\Tc) = 16 f_2 - 12 f_1^o + 9 f_0^o +
  h^{2, 2}_{3, 3} (\Tc) = f_0^+ - f_0^T + 2 f_0^b + 8 + h^{2,
  2}_{3, 3} (\Tc).
\end{equation}
If $\Tc$ is a hierarchical $(4,4)$-weighted subdivision, then by Theorem
\ref{thm:weighted}, we have  $h_{3, 3}^{2, 2}(\Tc)=0$. 


\begin{thebibliography}{10}

\bibitem{Berdinsky:2012}
Dmitry Berdinsky, Min-jae Oh, Tae-wan Kim, and Bernard Mourrain.
\newblock On the problem of instability in the dimension of a spline space over
  a T-mesh.
\newblock {\em Comput. Graph.}, 36(5):507--513, 2012.

\bibitem{b-htssg-88}
L.~Billera.
\newblock Homology theory of smooth splines: generic triangulations and a
  conjecture of {Strang}.
\newblock {\em Trans. Amer. Math. Soc.}, 310:325--340, 1988.

\bibitem{BuChSa10}
A.~Buffa, D.~Cho, and G.~Sangalli.
\newblock Linear independence of the {T}-spline blending functions associated
  with some particular {T}-meshes.
\newblock {\em Computer Methods in Applied Mechanics and Engineering},
  199:1437–1445, 2010.

\bibitem{CarEil99}
Henri Cartan and Samuel Eilenberg.
\newblock {\em Homological Algebra}.
\newblock Princeton University Press, 1999.

\bibitem{deBoor01}
Carl deBoor.
\newblock{\em A Practical Guide to Splines}.
\newblock Springer, 2001.

\bibitem{deng06}
Jiansong Deng, Falai Chen, and Yuyu Feng.
\newblock Dimensions of spline spaces over T-meshes.
\newblock {\em J. Comput. Appl. Math.}, 194(2):267--283, 2006.

\bibitem{Deng08TSPL22}
Jiansong Deng, Falai Chen, and Liangbing Jin.
\newblock Dimensions of biquadratic spline spaces over T-meshes.
\newblock arXiv:0804.2533v1.

\bibitem{deng08}
Jiansong Deng, Falai Chen, Xin Li, Changqi Hu, Weihua Tong, Zhouwang Yang, and
  Yuyu Feng.
\newblock Polynomial splines over hierarchical {T}-meshes.
\newblock {\em Graph. Models}, 70(4):76--86, 2008.

\bibitem{DoLyPe10}
Tor Dokken, Tom Lyche, and Kjell-Fredrik Pettersen.
\newblock Locally refined splines.
\newblock Preprint, 2010.

\bibitem{EhrRota93}
Richard Ehrenborg and Gian-Carlo Rota.
\newblock Apolarity and canonical forms for homogeneous polynomials.
\newblock {\em Eur. J. Comb.}, 14(3):157--181, 1993.


\bibitem{farin:book}
Gerald Farin. Curves and Surfaces for Computer Aided Geometric
Design: A Practical Guide, 5th Edition. Morgan Kaufmann, San Mateo, CA, 2001

\bibitem{Hatcher02}
Allen Hatcher.
\newblock {\em Algebraic Topology}.
\newblock Cambridge University Press, 2002.

\bibitem{HuDeFeChe06}
Zhangjin Huang, Jiansong Deng, Yuyu Feng, and Falai Chen.
\newblock New proof of dimension formula of spline spaces over {T}-meshes via
  smoothing cofactors.
\newblock {\em J. of Computational Mathematics}, 24(4):501–514, 2006.

\bibitem{hughes:CMAME2005}
Thomas J.R. Hughes, John A. Cottrell, Yuri Bazilevs.  Isogeometric analysis: CAD, finite elements, NURBS, exact geometry, 
and mesh refinement. Computer Methods in Applied Mechanics
and Engineering 194, 39-41, pp 4135-4195, 2005.

\bibitem{KungRota84}
Joseph P.~S. Kung and Gian-Carlo Rota.
\newblock The invariant theory of binary forms.
\newblock {\em Bull. Amer. Math. Soc. (N.S.)}, 10(1):27--85, 1984.

\bibitem{LiWaZha06}
Chong-Jun Li, Ren-Hong Wang, and Feng Zhang.
\newblock Improvement on the dimensions of spline spaces on {T}-mesh.
\newblock {\em Journal of Information \& Computational Science},
  3(2):235–244, 2006.

\bibitem{Li:2011:IDS:2038067.2038142}
Xin Li and Falai Chen.
\newblock On the instability in the dimension of splines spaces over t-meshes.
\newblock {\em Comput. Aided Geom. Des.}, 28(7):420--426, 2011.

\bibitem{Salmon1885}
G.~Salmon.
\newblock {\em Modern Higher Algebra}.
\newblock G.E.~Stechert and~Co., New York, 1885.
\newblock Reprinted~1924.

\bibitem{Schenck1997535}
Hal Schenck and Mike Stillman.
\newblock Local cohomology of bivariate splines.
\newblock {\em Journal of Pure and Applied Algebra}, 117-118:535 -- 548, 1997.

\bibitem{Sederberg04}
Thomas~W. Sederberg, David~L. Cardon, G.~Thomas Finnigan, Nicholas~S. North,
  Jianmin Zheng, and Tom Lyche.
\newblock {T}-spline simplification and local refinement.
\newblock {\em ACM Trans. Graph.}, 23(3):276--283, 2004.

\bibitem{Sederberg03}
Thomas~W. Sederberg, Jianmin Zheng, Almaz Bakenov, and Ahmad Nasri.
\newblock {T}-splines and {T}-nurccs.
\newblock {\em ACM Trans. Graph.}, 22(3):477--484, 2003.

\bibitem{Spanier66}
Edwin~H. Spanier.
\newblock {\em Algebraic Topology}.
\newblock MacGraw-Hill, 1966.

\bibitem{Wang01}
Ren-Hong Wang.
\newblock {\em Multivariate Spline Functions and Their Applications}.
\newblock Kluwer Academic Publishers, 2001.

\end{thebibliography}

\appendix
\section{Combinatorial properties}
We recall some well-known enumeration results for a T-mesh of an
axis-aligned rectangular domain $\Omega$.

\begin{lemma}
\label{lem:nbf}{\tmdummy}

\begin{itemizedot}
  \item $f_2 = f_0^+ + \frac{1}{2} f_0^T + \frac{1}{2} f_0^b - 1$
  
  \item $f_1^o = 2 f_0^+ + \frac{3}{2} f_0^T + \frac{1}{2} f_0^b - 2$
  
  \item $f_0^o = f_0^+ + f_0^T $
\end{itemizedot}
\end{lemma}

\begin{proof}
Each face $\sigma \in \Tc_2$ is a rectangle with 4 corners. If we
count these corners for all cells in $\Tc_2$, we enumerate 4 times
the crossing vertices, 2 times the $T$-vertices which are interior or on the
boundary and one time the corner vertices of $\Omega$. This yields the
relation
\[ 4 f_2 = 4 f_0^+ + 2 (f_0^T + (f_0^b - 4)) + 4. \]
Each interior edge $\tau \in \Tc_1^o$ has two end points.
Counting these end points for all interior edges, we count 4 times the
crossing vertices, 3 times the $T$-vertices which are interior and one time
the $T$-vertices on the boundary:
\[ 2 f_1^o = 4 f_0^+ + 3 f_0^T + (f_0^b - 4) . \]
Finally, as an interior vertex is a crossing vertex or a $T$-vertex, we have
\[ f_0^o = f_0^+ + f_0^T . \]
\end{proof}

\section{Complexes and homology\label{sec:homology}}

Let us recall here the basic properties that we will need on complexes of
vector spaces. Given a sequence of $\mathbb{K}$-vectors spaces $A_i$, $i = 0,
\ldots, l$ and linear maps $\partial_i : A_i \rightarrow A_{i + 1}$, we say
that we have a complex
\[ \mathcal{A}: A_l \rightarrow A_{l - 1} \rightarrow \cdots A_i \rightarrow
  A_{i - 1} \rightarrow \cdots A_1 \rightarrow A_0 \]
if $\tmop{im} \partial_i \subset \ker \partial_{i - 1}$.

\begin{definition}
 The $i^{\tmop{th}}$ homology $H_i (\mathcal{A})$ of $\mathcal{A}$ is $\ker
 \partial_{i - 1} / \tmop{im} \partial_i$ for $i = 1, \ldots, l$.
\end{definition}

The complex $\mathcal{A}$ is called {\tmem{exact}} (or an exact sequence) if
$H_i (\mathcal{A}) = 0$ (ie. $\tmop{im} \partial_i = \ker \partial_{i - 1}$)
for $i = 1, \ldots, l$. If the complex is exact and $A_l = A_0 = 0$, we have
\[ \sum_{i = 1}^{l - 1} (- 1)^i \dim A_i = 0. \]
Given complexes $\mathcal{A}= (A_i)_{i = 0, \ldots, l}$ $\mathcal{B}=
(\mathcal{B}_i)_{i = 0, \ldots, l}$, $\mathcal{C}= (C_i)_{i = 0, \ldots, l}$
and exact sequences
\[ 0 \rightarrow A_i \rightarrow B_i \rightarrow C_i \rightarrow 0 \]
for $i = 0, \ldots, l$, we have a long exact sequence
\cite{CarEil99}, \cite{Spanier66}[p. 182]:
\[ \cdots \rightarrow H_{i + 1} (\mathcal{C}) \rightarrow H_i (\mathcal{A})
  \rightarrow H_i (\mathcal{B}) \rightarrow H_i (\mathcal{C}) \rightarrow
  H_{i - 1} (\mathcal{A}) \rightarrow \cdots \]

\section{Dual topological complex}

The dual complex $\Tc^{\star}$ of the subdivision $\Tc$, is
such that we have the following properties.
\begin{itemizedot}
 \item a face $\sigma \in \Tc_2$ \ is a vertex of the dual complex
 $\Tc^{\star}$.
 
 \item an edge of $\Tc^{\star}$ is connecting two elements $\sigma,
 \sigma' \in \Tc_2$ if they share a common (interior) edge $\tau \in
 \Tc_1^o .$ Thus it is identified with the edge $\tau$ of
 $\Tc$ between $\sigma, \sigma'$;
 
 \item a face of $\Tc^{\star}$ corresponds to an elements $\gamma \in
 \Tc_0^o$. It is either a triangle if $\gamma$ is a $T$-junction or a
 quadrangle if $\gamma$ is a crossing vertex.
\end{itemizedot}
Notice that the boundary cells of $\Tc$ correspond to boundary
vertices of $\Tc^{\star}$. They are connected by boundary edges which
belong to a single face of $\Tc^{\star}$.

\section{Topological chain complex}
In this appendix section, we recall the main properties of the
topological chain complex 
{\small
\[ \begin{array}{lllllllll}
 \mathfrak{R}_{m, m'} (\Tc^o) : &  & \bigoplus_{\sigma \in
     \Tc_2} [\sigma] R_{m, m'} & \rightarrow & \bigoplus_{\tau \in
     \Tc_1^o} [\tau] R_{m, m'} & \rightarrow & \bigoplus_{\gamma \in
     \Tc_0^o} [\gamma] R_{m, m'} & \rightarrow & 0
   \end{array} \]
}
where 
\begin{itemize}
 \item $\forall \gamma \in \Tc_o^o$, $\partial_0 ([\gamma]) = 0$.
 \item  $\forall \tau=[\gamma_{1},\gamma_{2}] \in \Tc_1^o$, $\partial_1 
   ([\gamma_{1},\gamma_{2}]) = [\gamma_{1}]-[\gamma_{2}]$ with $[\gamma] \equiv 0$ iff
   $\gamma\in\partial \Omega$.
 and 
 \item $\forall \sigma \in \Tc_2^o$ with its counter-clockwise
  boundary formed by the edges $[\gamma_{1},\gamma_{2}], \ldots, [\gamma_{l},\gamma_{1}]$,
  $\partial_2 (\sigma) = [\gamma_{1},\gamma_{2}]  + \cdots +
  [\gamma_{l},\gamma_{1}]$ with $[\gamma,\gamma'] \equiv 0$ iff $\gamma,\gamma'\in\partial \Omega$.
\end{itemize}
We assume that $\Omega$ is simply connected and $\Omega^{o}$ is connected.
We prove that $\mathfrak{R}_{m, m'} (\Tc^o)$ is acyclic on a $T$-mesh of $\Omega$.

\begin{proposition}\label{proph0r}
  $H_0 (\mathfrak{R}_{m, m'} (\Tc^o)) = 0$.
\end{proposition}
\begin{proof}
  Let $\gamma \in \Tc_0^o$. There is a sequence of edges $\tau_0 =
  \gamma_0 \gamma_1$, \ $\tau_1 = \gamma_1 \gamma_2$, $\ldots$, \ $\tau_l =
  \gamma_l \gamma_{l + 1}$, such that $\tau_i \in \Tc_1^o$, $\gamma_0
  \not\in \Tc_0^o$ and $\gamma_{l + 1} = \gamma$. Then
  \[ \partial_1 ([\tau_0] + \cdots + [\tau_l]) = [\gamma_1] - [\gamma_0] +
     \cdots + [\gamma_{l + 1}] - [\gamma_l] = [\gamma] \]
  since $[\gamma_0] = 0$ and $[\gamma_{l + 1}] = [\gamma]$. Multiplying by any
  element in $R_{m, m'}$, we get that $[\gamma] R_{m, m'} \subset \tmop{im} \partial_1
  $ and thus $H_0 (\mathfrak{R}_{m, m'} (\Tc^o)) = 0$.
\end{proof}

\begin{proposition}\label{proph1r}
  $H_1 (\mathfrak{R}_{m, m'} (\Tc^o)) = 0$.
\end{proposition}

\begin{proof}
  Let $p = \sum_{\tau \in \Tc_1^o} [\tau] p_{\tau} \in \ker
  \partial_1$ with $p_{\tau} \in R_{m, m'}$. Let us prove that $p$ is in the
  image of $\partial_2$. For each $\gamma \in \Tc_0^o$ and each edge
  $\tau$ which contains $\gamma$, we have $\sum \varepsilon_{\tau} p_{\tau} =
  0$ with $\varepsilon_{\tau} = 1$ if $\tau$ ends at $\gamma$,
  $\varepsilon_{\tau} = - 1$ if $\tau$ starts at $\gamma$ and
  $\varepsilon_{\tau} = 0$ otherwise.
  
  For any $\sigma \in \Tc_2^o$ and $\tau \in \Tc_1^o$, we
  define $\varepsilon_{\sigma, \tau} = 1$ if $\tau$ is oriented
  counter-clockwise on the boundary of $\sigma$, $\varepsilon_{\sigma, \tau} =
  - 1$ if $\tau$ is oriented counter-clockwise on the boundary of $\sigma$ and
  $\varepsilon_{\sigma, \tau} = 0$ otherwise.
  
  For any oriented edge of the dual graph $\Tc^{\star}$ from $\sigma'$
  to $\sigma$, let us define \ $\partial_1^{\star} ([\sigma', \sigma]) =
  \varepsilon_{\sigma, \tau} p_{\tau}$. Notice that $\partial_1^{\star}
  ([\sigma, \sigma']) = \varepsilon_{\sigma', \tau} p_{\tau} = -
  \varepsilon_{\sigma, \tau} p_{\tau} = - \partial_1^{\star} ([\sigma',
  \sigma])$, since the orientation of $\tau$ on the boundary of $\sigma$ and
  $\sigma'$ are opposite.
  
  Let $\sigma_0$ be the 2-face of $\Tc$ with the lowest left corner
  for the lexicographic ordering. We
  order the cells $\sigma \in \Tc_2$ according to their distance to
  $\sigma_0$ in this dual graph $\Tc^{\star}$. We define an element $q
  = \sum_{\sigma \in \Tc_2} q_{\sigma} [\sigma]$ \ where $q_{\sigma}
  \in R_{m, m'}$ by induction using this order, as follows:
  \begin{itemizedot}
    \item $q_{\sigma_0} = 0$;
    
    \item For any $\sigma, \sigma' \in \Tc_2$, if $\sigma > \sigma'$
    and $\sigma$ and $\sigma'$ share a common edge $\tau$, then $q_{\sigma} =
    q_{\sigma'} + \partial_1^{\star} ([\sigma', \sigma])$.
  \end{itemizedot}
  Thus, if $[\sigma_0, \sigma_1], [\sigma_1, \sigma_2], \ldots, [\sigma_{k -
  1}, \sigma_k]$ is a path of $ \Tc^{\star}$ connecting $\sigma_0$ to
  $\sigma_k = \sigma$ with $\sigma_{i + 1} > \sigma_i$ then $ q_{\sigma} =
  \sum_{i = 0}^{k - 1} \partial_1^{\star} ([\sigma_i, \sigma_{i + 1}])$. Let
  us prove that this definition does not depend on the chosen path between
  $\sigma_0$ and $\sigma$.

  We first show that for any face {\tmem{$\gamma^{\star}$}} of \
  $\Tc^{\star}$ attached to a vertex $\gamma$, if its
  counter-clockwise boundary is formed by the edges $[\sigma, \sigma'],
  [\sigma', \sigma''], \ldots, [\sigma''', \sigma]$ corresponding to the edges
  $\tau, \tau', \tau'', \ldots$ of $\Tc$ containing $\gamma$, then
  \[ \partial_1^{\star} ([\sigma, \sigma']) + \partial_1^{\star} ([\sigma',
     \sigma''] + \cdots + \partial_1^{\star} ([\sigma''', \sigma]) =
     \varepsilon_{\sigma, \tau} p_{\tau} + \varepsilon_{\sigma', \tau'}
     p_{\tau'} + \varepsilon_{\sigma'', \tau''} p_{\tau''} + \cdots = 0. \]
  By changing the orientation of an edge $\tau$, we replace $p_{\tau}$ by $-
  p_{\tau}$ and \ $\varepsilon_{\sigma, \tau}$ by $- \varepsilon_{\sigma,
  \tau}$ so that the quantity $\varepsilon_{\sigma, \tau} p_{\tau}$ is not
  changed. Thus we can assume that all the edges $\tau, \tau', \tau'', \ldots$
  are pointing to $\gamma$. As $p \in \ker \partial_1$, we have $p_{\tau} +
  p_{\tau'} + p_{\tau''} + \cdots = 0$.
  
  Now as the cells $\sigma, \sigma', \sigma'', \ldots, \sigma$ are ordered
  counter-clockwise around $\gamma$ and as the edges are pointing to $\gamma$,
  we have $\varepsilon_{\sigma, \tau} = \varepsilon_{\sigma', \tau'} =
  \varepsilon_{\sigma'', \tau''} = \cdots = 1$, so that the sum
  $\partial_1^{\star} ([\sigma, \sigma']) + \partial_1^{\star} ([\sigma',
  \sigma''] + \cdots + \partial_1^{\star} ([\sigma''', \sigma])$ over a the
  boundary of a face {\tmem{$\gamma^{\star}$}} of \ $\Tc^{\star}$ is
  0.
  
  By composition, for any loop of $\Tc^{\star}$, the sum on the
  corresponding oriented edges is $0$. This shows that the definition of
  $q_{\sigma}$ does not depend on the oriented path from $\sigma_0$ to
  $\sigma$.
  
  By construction, we have
  \[ \partial_2 (q) = \sum_{\sigma \in \Tc_2} ( \sum_{\tau \in
     \Tc_1^o} \varepsilon_{\sigma, \tau} q_{\sigma} [\tau]) =
     \sum_{\tau \in \Tc_1^o} ( \sum_{\sigma \in \Tc_2}
     \varepsilon_{\sigma, \tau} q_{\sigma}) [\tau] . \]
  For each interior edge $\tau \in \Tc_1^o$, there are two faces
  $\sigma_1 > \sigma_2$, which are adjacent to $\tau .$ Thus, we have
  $\varepsilon_{\sigma_1, \tau} = - \varepsilon_{\sigma_2, \tau} $ and
  $q_{\sigma_1} = q_{\sigma_2} + \varepsilon_{\sigma_1, \tau} p_{\tau}$. We
  deduce that
  \[ ( \sum_{\sigma} \varepsilon_{\sigma, \tau} q_{\sigma}) =
     \varepsilon_{\sigma_1, \tau} q_{\sigma_1} + \varepsilon_{\sigma_2, \tau}
     q_{\sigma_2} = \varepsilon_{\sigma_1, \tau} (q_{\sigma_2} +
     \varepsilon_{\sigma_1, \tau} p_{\tau}) + \varepsilon_{\sigma_2, \tau}
     q_{\sigma_2} = p_{\tau} . \]
  This shows that $\partial_2 (q) = p$. In other words, $\tmop{im} \partial_2
  = \ker \partial_1$ and $H_1 (\mathfrak{R}_{m, m'} (\Tc^o)) = 0$.
\end{proof}

\begin{proposition}\label{proph2r} If $\Omega$ has one connected
  component. Then $H_2 (\mathfrak{R}_{m, m'} (\Tc^o)) = R_{m, m'}.$
\end{proposition}

\begin{proof}
  An element of $H_2 (\mathfrak{R}_{m, m'} (\Tc^o)) = \ker \partial_2$
  is a collection of polynomials $(p_{\sigma})_{\sigma \in \Tc_2}$ such
  that $p_{\sigma} \in R_{m, m'}$ and $p_{\sigma} = p_{\sigma'}$ if $\sigma$
  and $\sigma'$ share an (internal) edges. As $\Tc$ is a subdivision
  of a rectangle $D_0$, all faces $\sigma \in \Tc_2$ share pairwise an
  edge. Thus $p_{\sigma} = p_{\sigma'}$ for all $\sigma, \sigma' \in
  \Tc_2$ and $H_2 (\mathfrak{R}_{m, m'} (\Tc)) = R_{m, m'}$.
\end{proof}

Notice that by counting the dimensions in the exact sequence $\mathfrak{R}_{m,
m'} (\Tc)$, we recover the well-known Euler formula: $f_2 - f_1 + f_0
= 1$ (the domain $\Omega$ has one connected component).

\end{document}